\documentclass[11pt]{amsart}

\usepackage{fullpage}
\usepackage{amssymb}
\usepackage{amsmath}
\usepackage{amsxtra}
\usepackage{amscd}
\usepackage{graphicx}
\usepackage{amsfonts}
\usepackage{pb-diagram}
\usepackage{float}
\usepackage{enumerate}
\usepackage{dsfont}
\usepackage{multirow}
\usepackage{color}
\usepackage{rotating}
\usepackage{url}
\usepackage{enumitem}
\usepackage[bookmarks=false]{hyperref}
\definecolor{darkgreen}{rgb}{0,0.4,0.1}
\definecolor{darkpurple}{rgb}{0.7,0.0,0.4}
\hypersetup{
   colorlinks=true,       	
   linkcolor=darkpurple,     
   citecolor=darkgreen,    
   filecolor=magenta,      
   urlcolor=blue			
}

\numberwithin{equation}{section}
\newtheorem{thm}{Theorem}[section]

\newtheorem{cor}[thm]{Corollary}
\newtheorem{prop}[thm]{Proposition}

\newtheorem{conj}[thm]{Conjecture}

\newtheorem{rem}[thm]{Remark}

\def\N{{\mathbb N}}
\def\Z{{\mathbb Z}}

\def\C{{\mathbb C}}

\title{The BMM symmetrising trace conjecture for the 
exceptional 2-reflection groups of rank 2}
\author{Christina Boura}
\address{Laboratoire de Math\'ematiques UVSQ, B\^atiment Descartes, 45 avenue des \'Etats-Unis,  78035 Versailles cedex, France.}
\email{christina.boura@uvsq.fr}

\author{Eirini Chavli}
\address{Institut f\"ur Algebra und Zahlentheorie, Universit\"at Stuttgart, Pfaffenwaldring 57, 
70569 Stuttgart, Germany.}
\email{eirini.chavli@mathematik.uni-stuttgart.de}

\author{Maria Chlouveraki}
\address{Laboratoire de Math\'ematiques UVSQ, B\^atiment Fermat, 45 avenue des \'Etats-Unis,  78035 Versailles cedex, France.}
\email{maria.chlouveraki@uvsq.fr}


\thanks{We would like to thank the algebra school GABY in Milan,
where this project started. 
We also thank Ivan Marin for pointing us out further reasons for which the group $G_{13}$ is special, and Gunter Malle for checking this article.}
\begin{document}

\maketitle

\begin{center}
	\emph{Dedicated to Michel Brou\'e}
\end{center}

\begin{abstract}
	We prove the symmetrising trace conjecture of Brou\'e, Malle and Michel for the generic Hecke algebra associated to the exceptional  irreducible complex reflection group $G_{13}$. Our result completes the proof of the conjecture for the exceptional 2-reflection groups of rank 2.
	\end{abstract}
\section{Introduction}\label{Introduction}

Real reflection groups are finite groups of real matrices generated by reflections. They are also known as finite Coxeter groups. Complex reflection groups are finite groups of complex matrices  generated by pseudo-reflections, that is, non-trivial elements that fix a hyperplane pointwise. They provide thus a natural generalisation of real reflection groups.

The classification of complex reflection groups consists of classifying the irreducible complex reflection groups. These either belong to an infinite series $G(de,e,n)$, with $d,e,n \in \N^*$, or are one of the 34 exceptional groups $G_4,\,G_5,\,\ldots,\,G_{37}$. Among the latter, the first 19, that is, groups $G_4,\,G_5,\,\ldots,G_{22}$ are of rank $2$, \emph{i.e.,} they consist of $2 \times 2$ matrices. 
The quotient of each of these groups by its centre is the rotation group of a platonic solid.

Some complex reflection groups are generated by actual reflections, that is, pseudo-reflections of order $2$. These are called $2$-reflection groups. The $2$-reflection groups include all finite Coxeter groups and certain non-real groups of exceptional type: $G_{12}$, $G_{13}$, $G_{22}$, $G_{24}$, $G_{27}$, $G_{29}$, $G_{31}$, $G_{33}$ and $G_{34}$.

There are only three exceptional $2$-reflection groups of rank $2$, $G_{12}$, $G_{13}$ and $G_{22}$.
Most exceptional groups of rank $2$ are well-generated, that is, they can be generated by $2$ pseudo-reflections. However, this is not the case for the three groups in question, each of them needing at least $3$ generating reflections.
Among the $3$ groups, only $G_{13}$ has generators that do not belong to the same conjugacy class. In fact, it is the unique non-real exceptional $2$-reflection group with this property (the only real one being the Weyl group $F_4$).

There are many properties of real reflection groups that generalise to the complex ones, and many that we hope to generalise. One of the latter is that some of the complex reflection groups could play the role of Weyl groups of certain, as yet mysterious, objects that generalise finite reductive groups, the ``Spetses'' \cite{BMM, MalleS}. Inspired by this idea, Brou\'e, Malle and Rouquier \cite{BMR} associated to each complex reflection group two objects classically associated to real reflection groups: a braid group and a Hecke algebra. Since then the theories of braid groups and of Hecke algebras associated with complex reflection groups have become subjects of study in their own right, with many connections with other algebraic structures, such as Cherednik algebras or quantum groups, and other areas of mathematics, such as knot theory. 

When it comes to braid groups, $2$-reflection groups become important because, as it can be observed through a case-by-case analysis, any braid group associated with a complex reflection group is isomorphic to the braid group of a $2$-reflection group.
Another interesting fact, as we will see later in the paper, is that the braid group associated with $G_{13}$ is isomorphic to the braid group of the Weyl group $G_2$.

In this paper, we will focus on Hecke algebras associated with complex reflection groups. 
There are two fundamental conjectures about their structure, both known to hold for real reflection groups. Even without being proven, these conjectures have been 
assumed to hold 
in every paper where Hecke algebras have been used in the last two decades.
The first one is the Brou\'e--Malle--Rouquier  \cite{BMR}  freeness conjecture and states that the generic Hecke algebra $\mathcal{H}(W)$ associated to a complex reflection group $W$ is a free module over its ring of definition of rank equal to the order of $W$. This conjecture is now a theorem thanks to the work of many people.

The second one is the 
 Brou\'e--Malle--Michel \cite{BMM} symmetrising trace conjecture and states that there exists a symmetrising trace function on $\mathcal{H}(W)$ that satisfies certain canonicality conditions. The existence of such a trace gives us a lot of insight into the modular representation theory of the Hecke algebra, that is, when  its parameters specialise to complex numbers. Until recently, besides the finite Coxeter groups, this conjecture was known to hold for the exceptional groups $G_4,\,G_{12},\,G_{22}$ and $G_{24}$ \cite{MM10, Mar46}. Less than a year ago we proved this conjecture for the groups $G_4, \,G_5,\, G_6,\, G_7,\,G_8$ in \cite{BCCK}. A symmetrising trace exists also on the Hecke algebra associated with $G(de,e,n)$ \cite{BreMa, MM98}, but it is not yet known whether it satisfies all canonicality conditions.
 
 From now on, we will only discuss about exceptional  non-real reflection groups.
The techniques that we developed in \cite{BCCK} allow us to work with the Hecke algebra of any group whose basis has a ``nice'' form, 
such as the ones the second author has constructed for groups  $G_4, \,\ldots,\,G_{15}$ \cite{Ch17}; the basis is nice in two senses: it has an inductive nature and it involves powers of a central element $z$, which is the image in the Hecke algebra of the generator of the centre of the corresponding braid group. We expect the Hecke algebras of  all groups of rank $2$ to have such a basis, and
this is why we have restricted ourselves to the study of these groups. Now, as far as $2$-reflection groups are concerned, these have the following particularity: almost all of them have a generic Hecke algebra that depends only on $2$ parameters (and these $2$ parameters can be easily become one), a fact which  makes  calculations with these algebras much easier ---  in our work, the number of parameters is the main factor affecting computational difficulty. Of course, the order of the group, which is equal to the rank of the associated generic Hecke algebra, is also an important factor in the difficulty of calculations; this is why, for the moment, the BMM symmetrising trace conjecture has been so far proved for the smallest $2$-reflection groups and the smallest groups of rank $2$.

However, there is one exceptional $2$-reflection group that is an exception to the rule of a $2$-parameter generic Hecke algebra: this is the group $G_{13}$. Actually, the number of parameters of the generic Hecke algebra associated with a complex reflection group $W$ is equal to 
$\sum_{{\rm s} \in \mathcal{S}/ \sim} {\rm order}({\rm s})$, where $\mathcal{S}$ is a minimal generating set of reflections for $W$ and $ \sim$ stands for conjugacy in the group.
The generators of all other $2$-reflection groups belong to the same conjugacy class, but in $G_{13}$, we have $2$ distinct conjugacy classes of generators. Since all generators are of order $2$ (it is a $2$-reflection group), the generic Hecke algebra of $G_{13}$ depends on $4$ parameters. This difficulty was also the main obstacle for proving the BMR freeness conjecture with the methodology used for all other $2$-reflection groups by Marin and Pfeiffer \cite{MaPf}. 

In this article, we prove the BMM symmetrising trace conjecture for the group $G_{13}$, thus establishing its validity for all exceptional $2$-reflection groups  of rank $2$. Our methodology resembles the one we used in \cite{BCCK} for $G_4,\ldots,G_8$, but has some important computational differences, which are 
detailed in Section \ref{Proof}. In a few words, we use the basis, denoted by $\mathcal{B}_{13}$, for $\mathcal{H}(G_{13})$ of \cite{Ch17}. 
We define a linear map $\tau$ on  $\mathcal{H}(G_{13})$ by setting  $\tau(b):=\delta_{1b}$ for $b \in \mathcal{B}_{13}$ (we have $1 \in \mathcal{B}_{13}$). We show that $\tau$ is the canonical symmetrising trace on $\mathcal{H}(G_{13})$ by proving that $\tau$ is a symmetrising trace, that is, the matrix 
$A:=(\tau(bb'))_{b,b' \in  \mathcal{B}_{13}}$ is symmetric and invertible, and that $\tau$ satisfies the canonicality conditions. In order to do this, we created a  \texttt{C++}/SAGE combination program that allowed us, for any generator $g$ of $\mathcal{H}(G_{13})$ and any element $b$ of the basis $\mathcal{B}_{13}$,  to write $gb$ as a linear combination of elements in $\mathcal{B}_{13}$ (in \cite{BCCK}, we only used  \texttt{C++} for this part). Note that $\tau(gb)$ is the coefficient of $1$ in the aforementioned linear combination. We were then able to compute with SAGE  the whole matrix 
$A$ using the inductive nature of the basis  $\mathcal{B}_{13}$. Knowing the entries of the matrix $A$ also allowed us to check the canonicality conditions.

All programs described in this  paper (or information about them) can be found on the project's webpage \cite{Web}.

\section{Preliminaries}
\subsection{Complex reflection groups}
Let $V$ be a finite dimensional $\mathbb{C}$-vector space.  A {\em complex reflection group} is a finite subgroup of $\mathrm{GL}(V)$ generated by \emph{pseudo-reflections}, that is,
non-trivial elements of  $\mathrm{GL}(V)$ whose fixed points in $V$ form a hyperplane.

 We call the \emph{field of definition} of $W$, and denote by $K(W)$, the field generated by the traces on $V$ of all the elements of $W$. Benard \cite{Ben} and Bessis \cite{Bes1} have proved that $K(W)$ is a splitting field for all representations of  $W$. If $K(W) \subseteq \mathbb{R}$, then $W$ is a finite Coxeter group, and
if $K(W)=\mathbb{Q}$, then $W$ is a Weyl group.

A complex reflection group $W$ is called \emph{irreducible} if it acts irreducibly on $V$; in this case,  we define the \emph{rank} of $W$ to be
 the dimension of $V$. 
 Since every complex reflection group is a direct product of irreducible ones (see, for example,  \cite[Proposition 1.24]{lehrer}), we can restrict to the study of irreducible complex reflection groups. 
  The following theorem is due to  Shephard and Todd \cite{ShTo} and is known as the ``Shephard-Todd classification''.

\begin{thm}\label{ShToClas} Let $W \subset \mathrm{GL}(V)$ be an irreducible complex
reflection group. Then one of
the following assertions is true:
\begin{itemize}
  \item There exists a positive integer $n$ such that $(W,V) \cong (\mathfrak{S}_n, \C^{n-1})$.\smallbreak
  \item There exist positive integers $d,e,n$ with $de>1$ and $(de,e,n)\neq(2,2,2)$ such
  that $(W,V) \cong (G(de,e,n),\C^n)$, where $G(de,e,n)$ is the group of all 
  $n \times n$ monomial matrices whose non-zero entries are ${de}$-th roots of unity, while the product of all non-zero
  entries is a $d$-th root of unity. \smallbreak
  \item $(W,V)$ is isomorphic to one of the 34 exceptional groups
  $G_n$, with $n=4,\ldots,37$ (ordered with respect to increasing rank).
\end{itemize}
\end{thm}
  Among the irreducible complex reflection groups we encounter the irreducible finite Coxeter groups. More precisely, $G(1,1,n) \cong A_{n-1}$, $G(2,1,n) \cong B_n$,  $G(2,2,n) \cong D_n$, $G(m,m,2) \cong I_2(m)$, 
 	$G_{23} \cong H_3$,  $G_{28}  \cong  F_4$, $G_{30}  \cong H_4$, $G_{35}  \cong  E_6$, $G_{36}  \cong  E_7$, $G_{37}  \cong E_8$.

If $W$ is an irreducible complex reflection group of rank $r$, then $W$ is called \emph{well-generated} if it can be generated by $r$  pseudo-reflections. 
If this is not the case, then $W$ can be generated by $r+1$ pseudo-reflections. The well-generated groups are $G(d,1,r)$, $G(e,e,r)$ and all exceptional groups except for 
$G_7$, $G_{11}$, $G_{12}$, $G_{13}$, $G_{15}$, $G_{19}$, $G_{22}$ and $G_{31}$.

Now, in this paper, we will be interested in the intersection of the following two families of complex reflection groups:
\begin{itemize}
\item\textbf{The exceptional 2-reflection groups:}  These groups are generated by actual reflections,  that is, pseudo-reflections  of order $2$. Apart from the Coxeter groups we mentioned above, this family contains the groups $G_{12}$, $G_{13}$, $G_{22}$, $G_{24}$, $G_{27}$, $G_{29}$, $G_{31}$, $G_{33}$ and $G_{34}$.\smallbreak
\item\textbf{The exceptional groups of rank 2:} This family contains the groups $G_4,\dots, G_{22}$. The quotient of each group by its centre is the  rotation group of a platonic solid. As a result, the exceptional groups of rank 2 are divided into three smaller families, according to whether the aforementioned quotient is the tetrahedral group (which is the alternating group $\mathfrak{A}_4$), octahedral group (which is the symmetric  group $\mathfrak{S}_4$), or icosahedral group (which is the alternating group $\mathfrak{A}_5$). More precisely, we have the \emph{tetrahedral family}, which includes the groups $G_4,\dots, G_7$, the \emph{octahedral family}, which includes the groups $G_8,\dots, G_{15}$, and the \emph{icosahedral family}, which includes the rest.
\end{itemize}

 The groups belonging to both
 families are $G_{12}$, $G_{13}$ and $G_{22}$, which have the following Coxeter-like presentations: 
$$\begin{array}{lcl}
G_{12} &=& \langle\, {\rm s}, {\rm t}, {\rm u}\,|\,{\rm s} {\rm t}{\rm u}{\rm s}= {\rm t} {\rm u} {\rm s} {\rm t},\,
{\rm t} {\rm u} {\rm s} {\rm t}={\rm u}{\rm s} {\rm t}{\rm u},\, {\rm s}^2= {\rm t}^2={\rm u}^2=1\,\rangle,\smallbreak\smallbreak\smallbreak\\
G_{13}&=& \langle\, {\rm s}, {\rm t}, {\rm u} \,|\,
{\rm t} {\rm u} {\rm s} {\rm t}={\rm u}{\rm s} {\rm t}{\rm u},\,
{\rm s} {\rm t}{\rm u}{\rm s}{\rm t}=
{\rm u}{\rm s} {\rm t}{\rm u}{\rm s} ,\, {\rm s}^2= {\rm t}^2={\rm u}^2=1\,\rangle,\smallbreak\smallbreak\smallbreak\\
G_{22} &=&  \langle\, {\rm s}, {\rm t}, {\rm u}\,|\, 
{\rm s} {\rm t}{\rm u}{\rm s}{\rm t}=
{\rm u}{\rm s} {\rm t}{\rm u}{\rm s},\, 
{\rm s} {\rm t}{\rm u}{\rm s}{\rm t} = {\rm t} {\rm u} {\rm s} {\rm t}{\rm u},\,
{\rm s}^2= {\rm t}^2={\rm u}^2=1\,\rangle.
\end{array}
$$ 

We observe that, together with $G_{31}$, these are the only not well-generated exceptional $2$-reflection groups.
However, among them, $G_{13}$ is the only one where not all generators belong to the same conjugacy class: only ${\rm t}$ and ${\rm u}$ are conjugate, with ${\rm t} =  ({\rm u}{\rm s} {\rm t}) {\rm u} ({\rm u}{\rm s} {\rm t})^{-1}$. In fact, the only other exceptional $2$-reflection group with this property is the Weyl group $F_4$.

\subsection{Braid groups} Let $W\subset\mathrm{GL}(V)$ be a complex reflection group. Let $\mathcal{R}$ denote the set of pseudo-reflections of $W$ and let $\mathcal{A}:=\{{\rm ker}({\rm s}-{\rm id}_V)\;|\;{\rm s}\in \mathcal{R}\}$ be the set of reflecting hyperplanes of $W$. Set $X:=V\setminus \cup_{H\in \mathcal{A}}H$. 
As shown in \cite[\S2B]{BMR}, we can always  restrict to the case where
$W$ is ``essential'', meaning that $\cap_{H\in \mathcal{A}}H=\{0\}$. 
Steinberg  \cite[Corollary 1.6]{steinberg} proved that the action of $W$ on $X$ is free.
 Therefore, it defines a Galois covering $X\rightarrow X/W$, which gives rise to the following exact sequence for every $x\in X$:
$$1\rightarrow \pi_1(X,x)\rightarrow \pi_1(X/W, \underline{x})\rightarrow W\rightarrow 1,$$
where $\underline{x}$ denotes the image of $x$  under the canonical surjection $ X\rightarrow X/W$.
Brou\'e, Malle and Rouquier \cite{BMR} defined the {\em pure braid group} $P(W)$ and the 
{\em braid group} $B(W)$ of $W$ as $P(W):=\pi_1(X,x)$ and $B(W):=\pi_1(X/W, \underline{x})$. Moreover, they associated to every element of $\mathcal{R}$ homotopy classes in $B(W)$ that we call \emph{braided reflections}  (for more details, one may refer to \cite[\S 2]{BMR}).

We know by \cite[Theorem 0.1]{Bes2} that the braid group $B(W)$ admits an \emph{Artin-like} presentation. The generators of this presentation are braided reflections, whose number is the minimal number of pseudo-reflections that generate $W$. 
The relations are homogeneous relations (both sides of the same length) between 
positive
words in the generating elements. We call these relations \emph{braid relations}. 
For the cases that we are interested in, we have:
$$\begin{array}{lcl}
B(G_{12}) &=& \langle\, {\bf s}, {\bf t}, {\bf u}\,|\,{\bf s} {\bf t}{\bf u}{\bf s}= {\bf t} {\bf u} {\bf s} {\bf t},\,
{\bf t} {\bf u} {\bf s} {\bf t}={\bf u}{\bf s} {\bf t}{\bf u}\, \rangle,\smallbreak\smallbreak\smallbreak\\
B(G_{13})&=& \langle\, {\bf s}, {\bf t}, {\bf u} \,|\,
{\bf t} {\bf u} {\bf s} {\bf t}={\bf u}{\bf s} {\bf t}{\bf u},\,
{\bf s} {\bf t}{\bf u}{\bf s}{\bf t}=
{\bf u}{\bf s} {\bf t}{\bf u}{\bf s} \, \rangle,\smallbreak\smallbreak\smallbreak\\
B(G_{22}) &=&  \langle\, {\bf s}, {\bf t}, {\bf u}\,|\, 
{\bf s} {\bf t}{\bf u}{\bf s}{\bf t}=
{\bf u}{\bf s} {\bf t}{\bf u}{\bf s},\, 
{\bf s} {\bf t}{\bf u}{\bf s}{\bf t} = {\bf t} {\bf u} {\bf s} {\bf t}{\bf u}\,\rangle.
\end{array}
$$
Note that, again, $G_{13}$ is of particular interest. It stands out for being the only one among the exceptional $2$-reflection groups of rank $2$ with braid relations  of different length (4 and 5). 
A less obvious fact, shown by Bannai  \cite{Ban}, is that $B(G_{13})$ is isomorphic to the Artin group of type $I_2(6)$, that is, the braid group of the Weyl group $G_2$. 
We have
$$B(G_2) = \langle\, {\bf a}, {\bf b} \,|\,  {\bf a} {\bf b} {\bf a} {\bf b} {\bf a} {\bf b}=
 {\bf b} {\bf a} {\bf b} {\bf a} {\bf b} {\bf a} \,\rangle.
$$
and an isomorphism is given by
${\bf a} \mapsto {\bf u} {\bf s},\, {\bf b} \mapsto {\bf u}{\bf s}{\bf t}({\bf u}{\bf s})^{-1}$,
with inverse ${\bf s} \mapsto {\bf b}{\bf a}  ({\bf b} {\bf a} {\bf b})^{-1}{\bf a},\,
 {\bf t} \mapsto  {\bf a}^{-1} {\bf b} {\bf a},\,  {\bf u} \mapsto ({\bf a} {\bf b} {\bf a})^{-1}
 {\bf b} ({\bf a} {\bf b} {\bf a})$
(the maps are taken from \cite[Proof of Proposition 7.4]{MarKramm}, but we have corrected the image of ${\bf s}$).

We assume now that $W$ is irreducible and we denote by $Z(W)$ the centre of $W$. The irreducibility of $W$ implies that $Z(W)$ is in bijection with a (finite) subgroup of $\C^{\times}$, and it is thus a cyclic group. Let $x$ be some fixed basepoint of $X$. We denote by $\boldsymbol{\pi}$ and $\boldsymbol{\beta}$ the homotopy classes of the loops $t\mapsto x \,{\rm exp}(2\pi i t)$ and $t\mapsto  x\,{\rm exp}(2\pi i t/|Z(W)|)$ respectively. Brou\'e, Malle and Rouquier \cite[Lemma 2.22 (2)]{BMR} proved that $\boldsymbol{\beta}$ and $\boldsymbol{\pi}$ belong to the centre of $B(W)$ and $P(W)$ respectively. Moreover,  they conjectured the following:
\begin{itemize}
	\item[(i)] The centre of $B(W)$ is cyclic and generated by $\boldsymbol{\beta}$ (proved by Bessis \cite[Theorem 12.8]{Bes3}).\smallbreak
	\item[(ii)] The centre of $P(W)$ is cyclic and generated by $\boldsymbol{\pi}$ (proved by Digne, Marin and Michel \cite[Theorem 1.2]{DMM}).
\end{itemize}

\subsection{Hecke algebras} 
\label{hecke}

A pseudo-reflection ${\rm s} \in \mathcal{R}$ is called \emph{distinguished} if its only nontrivial eigenvalue on $V$ equals ${\rm exp}(2\pi i /e_{\rm s})$, where
$e_{\rm s}$ denotes the order of ${\rm s}$ in $W$. Let $\mathcal{S}$ denote the set of the distinguished pseudo-reflections of $W$. For each ${\rm s} \in \mathcal{S}$ we choose a set of $e_{\rm s}$ indeterminates ${\rm v}_{{\rm s},1},\dots, {\rm v}_{{\rm s},e_{\rm s}}$, such that ${\rm v}_{{\rm s},k}={\rm v}_{{\rm t},k}$ if ${\rm s}$ and ${\rm t}$ are conjugate in $W$. We denote by $R(W)$ the Laurent polynomial ring $\Z[({\rm v}_{{\rm s},k},{\rm v}_{{\rm s},k}^{-1})_{{\rm s} \in \mathcal{S}}]$. The \emph{generic Hecke algebra} $\mathcal{H}(W)$ associated to $W$ with parameters $({\rm v}_{{\rm s},1},\dots, {\rm v}_{{\rm s},e_{\rm s}})_{{\rm s} \in \mathcal{S}}$ is the quotient of the group algebra $R(W)B(W)$  by the ideal generated by the elements of the form 
\begin{equation}
({\bf s}-{\rm v}_{{\rm s},1})({\bf s}-{\rm v}_{{\rm s},2})\dots ({\bf s}-{\rm v}_{{\rm s},e_{\rm s}}),
\label{Hecker}
\end{equation}
where ${\rm s}$ runs over 
$\mathcal{S}$ and ${\bf s}$ runs over the set of braided reflections associated to the pseudo-reflection ${\rm s}$. 
One may notice that it is enough to choose one relation of the form described in (\ref{Hecker}) per conjugacy class of $\mathcal{S}$, since  the corresponding braided reflections are conjugate in $B(W)$. 

We obtain an equivalent definition of $\mathcal{H}(W)$ if we expand the relations (\ref{Hecker}). 
More precisely, $\mathcal{H}(W)$ is the quotient of the group algebra $R(W)B(W)$  by the ideal generated by the elements of the form 
\begin{equation}
{\bf s}^{e_{\rm s}}-a_{{\rm s},e_{\rm s}-1}{\bf s}^{e_{\rm s}-1}-a_{{\rm s},e_{\rm s}-2}{\bf s}^{e_{\rm s}-2}-\dots-a_{{\rm s},0},
\label{Hecker2}
\end{equation}
where $a_{{\rm s},e_{\rm s}-k}:=(-1)^{k-1}f_k({\rm v}_{{\rm s},1},\dots,{\rm v}_{{\rm s},e_{\rm s}})$ with $f_k$ denoting the $k$-th elementary symmetric polynomial, for $k=1,\ldots,e_{\rm s}$. 
Therefore, in the presentation of $\mathcal{H}(W)$, we have two kinds of relations:
 the \emph{braid relations}, coming from the Artin-like presentation of $B(W)$, and the \emph{positive Hecke relations}: 
 \begin{equation}
{\bf s}^{e_{\rm s}}=a_{{\rm s},e_{\rm s}-1}{\bf s}^{e_{\rm s}-1}+a_{{\rm s},e_{\rm s}-2}{\bf s}^{e_{\rm s}-2}+\dots+a_{{\rm s},0}.
 \end{equation}
  We notice now that $a_{{\rm s},0}= (-1)^{e_{\rm s}-1}{\rm v}_{{\rm s},1}{\rm v}_{{\rm s},2}\dots {\rm v}_{{\rm s},e_{\rm s}} \in R(W)^{\times}$. Hence,  ${\bf s}$ is invertible in $\mathcal{H}(W)$ with
   \begin{equation}\label{invhecke}
 {\bf s}^{-1}=a_{s,0}^{-1}\,{\bf s}^{e_{\rm s}-1}-a_{s,0}^{-1}\,a_{{\rm s},e_{\rm s}-1}{\bf s}^{e_{\rm s}-2}-a_{s,0}^{-1}\,a_{{\rm s},e_{\rm s}-2}{\bf s}^{e_{\rm s}-3}-\dots- a_{s,0}^{-1}\,a_{{\rm s},1}.
 \end{equation}
 We call  relations \eqref{invhecke} the \emph{inverse Hecke relations}.
 
 Let us now see what happens in the case of the exceptional $2$-reflection groups of rank $2$. All generators of $G_{12}$ and $G_{22}$ belong to the same conjugacy class, and so the generic Hecke algebras associated to these groups are defined over a Laurent polynomial ring in two indeterminates, which we denote by ${\rm v}_1$ and ${\rm v}_2$. We have: 
  $$\mathcal{H}(G_{12}) =\left \langle\, {s}, {t}, {u}\,|\,{s} {t}{u}{s}= {t} {u} {s} {t},\,
{t} {u} {s} {t}={u}{s} {t}{u},\, \prod_{k=1}^2(s-{\rm v}_k)=\prod_{k=1}^2(t-{\rm v}_k)=
\prod_{k=1}^2(u-{\rm v}_k)=0
\,\right\rangle,$$
$$\phantom{aaa}\mathcal{H}(G_{22}) = \left\langle\, {s}, {t}, {u}\,|\, 
{s} {t}{u}{s}{t}=
{u}{s} {t}{u}{s},\, 
{s} {t}{u}{s}{t} = {t} {u} {s} {t}{u},\,
\prod_{k=1}^2(s-{\rm v}_k)=\prod_{k=1}^2(t-{\rm v}_k)=
\prod_{k=1}^2(u-{\rm v}_k)=0\,\right\rangle.
$$

However,
as we have already mentioned before, the generators of $G_{13}$ belong
 to two different conjugacy classes, so its generic Hecke algebra is defined over the Laurent polynomial ring 
$R(G_{13})=\Z[{\rm v}_{{\rm s},k},{\rm v}_{{\rm s},k}^{-1}, {\rm v}_{{\rm t},k},{\rm v}_{{\rm t},k}^{-1}]$ and has the following presentation:
  $$\mathcal{H}(G_{13}) =\left \langle\, {s}, {t}, {u}\,|\,
{t} {u} {s} {t}={u}{s} {t}{u},\,{s} {t}{u}{s}{t}=
{u}{s} {t}{u}{s},\,  \prod_{k=1}^2(s-{\rm v}_{{\rm s},k})=\prod_{k=1}^2(t-{\rm v}_{{\rm t},k})=
\prod_{k=1}^2(u-{\rm v}_{{\rm t},k})=0
\,\right\rangle.$$
Setting $a:={\rm v}_{{\rm s},1}+{\rm v}_{{\rm s},2}$,
$b:=-{\rm v}_{{\rm s},1}{\rm v}_{{\rm s},2}$,
$c:={\rm v}_{{\rm t},1}+{\rm v}_{{\rm t},2}$ and
$d:=-{\rm v}_{{\rm t},1}{\rm v}_{{\rm t},2}$, we can rewrite the presentation of 
 $\mathcal{H}(G_{13})$ as follows:
  \begin{equation}\label{g13}
  \mathcal{H}(G_{13}) =\left \langle\, {s}, {t}, {u}\,|\,
{t} {u} {s} {t}={u}{s} {t}{u},\,{s} {t}{u}{s}{t}=
{u}{s} {t}{u}{s},\,  s^2=as+b,\,t^2=ct+d,\,u^2=cu+d
\,\right\rangle.
\end{equation}
Note that only $b$ and $d$ are units in $R(G_{13})$. The inverse Hecke relations in this case are:
 \begin{equation}\label{invg13}
 s^{-1} = b^{-1}s-b^{-1}a,\;\;\;t^{-1} = d^{-1}t-d^{-1}c,\;\;\;u^{-1} = d^{-1}u-d^{-1}c.
 \end{equation}

\section{Conjectures}\label{Conjectures}
Let $W\subset GL(V)$ be a complex reflection group and let $\mathcal{H}(W)$ be its associated generic Hecke algebra over the Laurent polynomial ring $R(W)$, as defined in the previous section. 

\subsection{The BMR freeness conjecture}
In  \cite[\S 4]{BMR},
Brou\'e, Malle and Rouquier stated the following conjecture on the structure of $\mathcal{H}(W)$:

\begin{conj}\label{BMR free}
	The algebra $\mathcal{H}(W)$ is a free $R(W)$-module of rank  $|W|$.
	\end{conj}
	
Note that, thanks to the following result, which can be found in 	\cite[Proof of Theorem 4.24]{BMR} (for another detailed proof, one may also see \cite[Proposition 2.4]{Mar43}), 
	in order to prove the BMR freeness conjecture, it is enough to find a spanning set of $\mathcal{H}(W)$ over $R(W)$ consisting of $|W|$ elements.
	\begin{thm}\label{mention}
	 If $\mathcal{H}(W)$ is generated as $R(W)$-module by  $|W|$ elements, then it is a free 
	 $R(W)$-module of rank  $|W|$.
	\end{thm}
	
The BMR freeness conjecture is now a theorem.	
When stated, this conjecture was already known to hold for real reflection groups \cite[IV, \S 2]{Bou05} and for the complex reflection groups of the infinite series $G(de,e,n)$ \cite{ArKo, BM, Ar}. As far as the exceptional non-real reflection groups are concerned, we have the validity of the conjecture for
	\begin{itemize}
	\item the group $G_4$ by \cite{BM, Fun, Mar41, Ch18} (4 independent proofs); \smallbreak
\item the group $G_{12}$ by  \cite{MaPf}; \smallbreak
\item the groups $G_5,\ldots, G_{16}$ by \cite{Ch17, Ch18}; \smallbreak
\item the groups $G_{17},\,G_{18},\,G_{19}$ by \cite{Tsu} (with a computer method applicable to all rank $2$ groups);\smallbreak
\item the groups $G_{20}, \,G_{21}$ by \cite{MarNew};\smallbreak
\item the groups $G_{22},\ldots, G_{37}$ by \cite{Mar41, Mar43, MaPf}. 
	\end{itemize}

\subsection{The BMM symmetrising trace conjecture} 
In \cite[2.1]{BMM}, Brou\'e, Malle and Michel stated a second fundamental conjecture on the structure of $\mathcal{H}(W)$:
	\begin{conj}\label{BMM sym}
	There exists a linear map $\tau: \mathcal{H}(W)\rightarrow R(W)$ such that:
	\begin{enumerate}
		\item $\tau$ is a symmetrising trace, that is, we have $\tau(h_1h_2)=\tau(h_2h_1)$ for all  $h_1,h_2\in \mathcal{H}(W)$, and the bilinear map $\mathcal{H}(W)\times \mathcal{H}(W)\rightarrow R(W)$,  $(h_1,h_2)\mapsto\tau(h_1h_2)$ is non-degenerate. \smallbreak
	\item $\tau$ becomes the canonical symmetrising trace on $K(W)W$ (given by 
	$\tau(w)=\delta_{1w}$) when ${\rm v}_{{\rm s},k}$ specialises to ${\rm exp}(2\pi i k/e_{\rm s})$ for all ${\rm s} \in \mathcal{S}$ and $k=1,\ldots, e_{\rm s}$. \smallbreak
		\item  $\tau$ satisfies
	$$
	\tau(T_{b^{-1}})^* =\frac{\tau(T_{b\boldsymbol{\pi}})}{\tau(T_{\boldsymbol{\pi}})}, \quad \text{ for all } b \in B(W),
	$$
	where $b\mapsto T_{b}$ denotes the restriction of the natural surjection $R(W)B(W) \rightarrow \mathcal{H}(W)$ to $B(W)$ and
	$x \mapsto x^*$ the automorphism of  $R(W)$ given by ${\rm v}_{{\rm s},k} \mapsto {\rm v}_{{\rm s},k}^{-1}$.
	\smallbreak
	\end{enumerate}
	\end{conj}
	Note that, by  \cite[2.1]{BMM}, since the BMR freeness conjecture holds,  there exists at most
	one symmetrising trace satisfying Conditions (2) and (3) of Conjecture \ref{BMM sym}, meaning that $\tau$ is unique. We call $\tau$ the \emph{canonical symmetrising trace on  $\mathcal{H}(W)$}.

The BMM symmetrising trace conjecture holds for all real reflection groups \cite[IV, \S 2]{Bou05}  and	
for the following exceptional complex reflection groups:
\begin{itemize}
	\item the group $G_4$ by \cite{MM10, Mar46, BCCK} (3 independent proofs); \smallbreak
	\item the groups $G_5, G_6, G_7, G_8$ by \cite{BCCK}; \smallbreak
	\item the groups $G_{12}, G_{22}, G_{24}$ by \cite{MM10}.
\end{itemize}	
A map that satisfies Conditions (1) and (2) of the conjecture exists also for the groups of the infinite series $G(de,e,n)$ by \cite{BreMa, MM98};  it is still unclear though whether this symmetrising trace satisfies Condition (3).

In all the cases above, 
a suitable basis $\mathcal{B}(W)$ for $\mathcal{H}(W)$ was considered, so that:
\begin{enumerate}
\item[(i)] $1 \in \mathcal{B}(W)$; \smallbreak
\item[(ii)] each element of $\mathcal{B}(W)$ corresponds to a distinct element of $W$ when ${\rm v}_{{\rm s},k}$ specialises to ${\rm exp}(2\pi i k/e_{\rm s})$
for all ${\rm s} \in \mathcal{S}$ and $k=1,\ldots, e_{\rm s}$; \smallbreak
\item[(iii)]  $\tau(b)=\delta_{1b}$ for all $b \in \mathcal{B}(W)$.
\end{enumerate}	
This way, 	Condition (2) of Conjecture \ref{BMM sym} is satisfied, and only Conditions (1) and (3) have to be checked. In fact, Malle and Michel conjectured that there is always a subset of $\mathcal{H}(W)$ (not necessarily a basis) that satisfies properties (i)--(iii) \cite[Conjecture 2.6]{MM10}:
	
	\begin{conj}\textbf{``The lifting conjecture''}
	There exists a section $W \rightarrow \boldsymbol{W} \subset B(W)$, $w \mapsto \boldsymbol{w}$ of $W$ in $B(W)$ such that $1 \in \boldsymbol{W}$, and such that for any $\boldsymbol{w} \in \boldsymbol{W}$ we have $\tau(T_{\boldsymbol{w}}) = \delta_{1\boldsymbol{w}}$.
	\end{conj}
	
	If the above conjecture holds, then Condition (2) of Conjecture \ref{BMM sym} is obviously satisfied. If further the elements $\{T_{\boldsymbol{w}}\,|\,\boldsymbol{w} \in \boldsymbol{W}\}$ form an
	$R(W)$-basis of $\mathcal{H}(W)$, then, by \cite[Proposition 2.7]{MM10}, Condition (3) of Conjecture \ref{BMM sym} is equivalent to:
	\begin{equation}\label{extra con2}
	\tau(T_{\boldsymbol{w}^{-1}\boldsymbol{\pi}})=0, \quad \text{ for all } \boldsymbol{w} \in \boldsymbol{W} \setminus \{1\}.
	\end{equation}

	\section{Results}\label{Proof}
	
	The main result of this paper is the following:
	
	\begin{thm}
		The BMM symmetrising trace conjecture holds for $G_{13}$.	
		\end{thm}
	
Since the validity of the conjecture is already known for the groups $G_{12}$ and $G_{22}$ by \cite{MM10}, we conclude the following:	
	
	\begin{cor}
		The BMM symmetrising trace conjecture holds for the exceptional 2-reflection groups of rank 2.
	\end{cor}

The main difficulty of $G_{13}$ in comparison with $G_{12}$ and $G_{22}$ is that  not all its generators belong to the same conjugacy class. In fact, it is the only non-real exceptional $2$-reflection group that has this property. The presence of two conjugacy classes of generators was also the main obstacle for proving the BMR freeness conjecture for $G_{13}$ with the  methodology provided by Marin and Pfeiffer \cite{MaPf}, who proved the BMR freeness conjecture for all the other non-real exceptional $2$-reflection groups.

Recall now the presentation of the generic Hecke algebra $\mathcal{H}(G_{13})$ given by \eqref{g13}. Note that, since the generators of $G_{13}$ are of order $2$ and belong to two different conjugacy classes, $\mathcal{H}(G_{13})$ depends on four parameters $a$, $b$, $c$ and $d$; among them, only $b$ and $d$ are units in $R(G_{13})$.
In \cite{Ch17}, the second author proved the BMR freeness conjecture for  $G_{13}$ by providing the following explicit basis for  $\mathcal{H}(G_{13})$:
$$\mathcal{B}_{13}: =
\left\{
\begin{matrix}
\begin{array}{c|l}
z^kt^l, z^kut^l, z^kst^l, z^ktst^l, z^ksut^l, z^kust^l,\\  z^ktut^l, z^ktsut^l, z^ktust^l, z^kstst^l, z^kstut^l, z^kutst^l
\end{array}
\end{matrix}\,\,\,k=0,1,2,3,\,\,\,l=0,1\right\},
$$
where $z=T_{\boldsymbol{\beta}}$, the image of the generator $\boldsymbol{\beta}$ of the centre of $B(G_{13})$  inside $\mathcal{H}(G_{13})$. We know  by \cite[Appendix 2, Table 3]{BMR} that $z=(stu)^3=(tus)^3=(ust)^3$. We have $|\mathcal{B}_{13}|=|G_{13}|=96$ and $|Z(G_{13})|=4$.

Notice here that $1\in \mathcal{B}_{13}$ and that there exists a section $G_{13} \rightarrow \boldsymbol{G}_{13}$, $w \mapsto \boldsymbol{w}$ of $G_{13}$ in $B(G_{13})$ such that  $\mathcal{B}_{13} = \{T_{\boldsymbol{w}}\,|\,\boldsymbol{w} \in \boldsymbol{G}_{13}\}$.
Let $\tau: \mathcal{H}(G_{13})\rightarrow R(G_{13})$ be the $R(G_{13})$-linear map defined by $\tau(b)= \delta_{1b}$, for every $b\in \mathcal{B}_{13}$. Obviously, $\tau$    satisfies Condition (2) of Conjecture \ref{BMM sym}. We will  prove that it also satisfies Conditions (1) and (3).

First though, we will introduce some notation. 
For every $k \in \{0,1,2,3\}$ and $l \in \{0,1\}$, we set    
$$\mathcal{E}_{13}^{k,l}:  =\left\{
z^kt^l, z^kut^l, z^kst^l, z^ktst^l, z^ksut^l, z^kust^l,  z^ktut^l, z^ktsut^l, z^ktust^l, z^kstst^l, z^kstut^l, z^kutst^l
\right\}.$$
The sets $\mathcal{E}_{13}^{k,l}$ form a partition of $\mathcal{B}_{13}$. Moreover, we have $\mathcal{E}_{13}^{k,1}=\mathcal{E}_{13}^{k,0}t$ for all $k \in \{0,1,2,3\}$.
We can now order the basis $\mathcal{B}_{13}$, using the set $\mathcal{E}_{13}^{0,0}=\{1, u, s, ts, su, us,  tu, tsu, tus, sts, stu, uts\}$. 
We notice that every $i\in\{1,\dots, 96\}$ can be written
uniquely as $24k+12l+m$, where $k\in\{0,1,2,3\}$, $l \in \{0,1\}$ and $m\in\{1,\ldots, 12\}$. We write $\mathcal{B}_{13}=\{b_1,\dots,b_{96}\}$, where $b_i:=z^kb_mt^l$ for $b_m \in \mathcal{E}_{13}^{0,0}$ in the order written above. Note that we also have $b_i=z^kb_{12l+m}$. For example, $b_{72}=b_{24\cdot2+12 \cdot 1+12}=z^2b_{12}t=z^2utst = z^2b_{24}$.

\subsection{Proof of the first condition} In order to prove that $\tau$ is a symmetrising trace, we need to show that the matrix $A:=(\tau(b_ib_j))_{1 \leq i,j \leq 96}$ is symmetric and that its determinant is a unit in $R(G_{13})$. 
In this subsection, we will present the algorithm that we used to calculate the entries of the matrix $A$. We then checked that $A$ is symmetric and computed its determinant to be equal to $b^{512}d^{1032}$.
Our algorithm follows the same idea as the algorithm used in \cite{BCCK} for groups $G_5$, $G_6$, $G_7$ and $G_8$, but some important modifications were required. 
All programs described in this section (or information about them) can be found on the project's webpage \cite{Web}.

 By definition of $\tau$, for any $h \in \mathcal{H}(G_{13})$, $\tau(h)$ is the coefficient of $1$ when $h$ is expressed as a linear combination of the elements of $\mathcal{B}_{13}$.  From now on, we set $\mathfrak{R}:=\mathbb{Z}[a,b^{\pm1},c,d^{\pm 1}]$.
 
 \subsubsection{First step: The products $sb_j$, $tb_j$ and $ub_j$}\label{first step}
 Our first step
 is to express $b_ib_j$ 
as an $\mathfrak{R}$-linear combination of the elements  of $\mathcal{B}_{13}$ for $i \in \{2,3,13\}$ (that is, $b_i \in \{u,s,t\}$) and  $j\in\{1,\dots,96\}$.
Due to the nice form of $\mathcal{B}_{13}$, one only needs to compute these   linear combinations  for  $b_j\in \mathcal{E}^{k,0}_{13}$, $k=0,1,2,3$. This follows from the fact that every $b_j\in \mathcal{E}^{k,1}_{13}$ can be written as $b_j=b_{j-12}t$ with 
$b_{j-12} \in \mathcal{E}^{k,0}_{13}$.
Hence, if $b_ib_{j-12} = \sum_{n} \lambda_n b_n$ for some $\lambda_n \in \mathfrak{R}$, then $b_ib_j = \sum_{n} \lambda_n b_n t$. If now $b_n \in \mathcal{E}^{k',0}_{13}$ for some $k' \in \{0,1,2,3\}$, then 
 $b_n t=b_{n+12}$.   
 On the other hand, if  $n\in\mathcal{E}^{k',1}_{13}$ for some $k' \in \{0,1,2,3\}$, then $b_n t=b_{n-12}t^2=cb_{n-12}t+db_{n-12}=cb_{n}+db_{n-12}.$
This argument is taken into account in all the ``special cases'' mentioned below.

Let $b_i \in\{u,s,t\}$ and $b_j\in \mathcal{E}^{k,0}_{13}$ for some $k \in \{0,1,2,3\}$. When we try to express
$b_ib_j$ as an $\mathfrak{R}$-linear combination of the elements  of $\mathcal{B}_{13}$, we encounter the following cases:
\begin{itemize}
	\item $b_ib_j\in \mathcal{B}_{13}$. For example, $sb_{4}=sts=b_{10}$. \smallbreak
	\item $b_ib_j$ is an $\mathfrak{R}$-linear combination of two elements of  $\mathcal{B}_{13}$. We obtain these elements by applying one of the positive Hecke relations. For example,
	$tb_{4}=t^2s=cts+ds=cb_{4}+db_{3}$. \smallbreak
	\item $b_ib_j$ is one of the 25 ``special cases''. These are some equalities computed by hand which express a given $b_ib_j$  as a sum of other elements in $\mathcal{H}(G_{13})$. 
	In order not to ruin the text's coherence, we will discuss 
	a bit further about the special cases and give examples 
	of how we use them in the Appendix of this paper. 
	One can find the complete list of special cases for $G_{13}$
	on the project's webpage \cite{Web}.
		 \smallbreak
	\item Combinations of the above.
\end{itemize}	
	
	It is entirely possible to do all these calculations by hand --- examples are provided in the Appendix.
	Actually, the second author has done all these calculations by hand in order to prove the BMR freeness conjecture for $G_{13}$ in \cite{Ch17}, without however keeping track of the explicit coefficients. Since the calculations are time-consuming and mistakes can be made when computing the coefficients, 
  we created a computer  program in the language \texttt{C++} whose purpose was to produce the desired linear combination for every product $b_ib_j$, with $b_i \in \{u,s,t\}$, using the methodology described below and illustrated in the examples 
  found in the Appendix.
  The  \texttt{C++} program is similar to the one constructed for \cite{BCCK}, but it has the following important differences:
\begin{itemize}
\item[D1.] We  expand the basis $\mathcal{B}_{13}$ to a spanning set 
$\mathcal{C}_{13}$  
which includes the following four extra elements: 
$b_{97}:=tstsu$, 
$b_{98}:=tstsut$, 
$b_{99}:=zututs$, and 
$b_{100}:=zututst$. These elements appear only when we try to express 
$sb_j$ as an $\mathfrak{R}$-linear combination of the elements of the basis (and, as it turns out, only for the products $sb_8$, $sb_{12}$, $sb_{20}$ and $sb_{24}$, although  knowing this in advance is not important).
 If we wanted the \texttt{C++} program to break further these elements, then we would need to introduce many more ``special cases''. To avoid this, we ask the program to
 simply express $sb_j$ as a linear combination of the  elements of  $\mathcal{C}_{13}$. We will see in a while how we continue from there. \smallbreak
 
 \item[D2.] We do not ask the program to compute the products $sb_{81}$ and $sb_{93}$, because again it would be computationally  complicated, while these cases can be treated separately.\smallbreak

\item[D3.] We know that the algorithm is not heuristic, since we can see, from the computations by hand, that every special case eventually leads to a linear combination of elements of $\mathcal{C}_{13}$.\smallbreak

\item[D4.] We do not need the braid relations as input, because we applied them whenever necessary by hand when we produced the ``special cases''.

\end{itemize}  
    
   Let us now describe the \texttt{C++} program.
  The inputs of the  program are the following:  
 \begin{itemize}
\item[I1.] The set $\mathcal{C}_{13}= \mathcal{B}_{13} \cup \{b_{97}, b_{98}, b_{99}, b_{100}\}$.\smallbreak
\item[I2.] The positive and inverse Hecke relations. \smallbreak  
\item[I3.] The 25 special cases.
\end{itemize}

To start with, we define a list $L$ of  elements of the form $\lambda \,T_b$, where $\lambda \in \mathfrak{R}$ and $b \in B(G_{13})$. We  initialise $L$ with the element $b_ib_j$, \emph{i.e.,} $L = [b_ib_j]$, with $i \in \{2,3,13\}$ and $j \in \{1,\ldots,96\}$ (we exclude the cases $(i,j)=(3,81)$ and $(i,j)=(3,93)$). The algorithm is iterative. At each iteration, the initial element is  expanded and replaced by its summands.  For example, suppose that the input is $b_{13}b_{13} = t^2$. So in the beginning $L = [t^2]$. During the first iteration, by the positive Hecke relation, $t^2$ is equal to $ct + d$, so the list $L$ is updated to become $L = [ct,d]$.

The same procedure is repeated at each iteration until all the elements of $L$ are of the form $\lambda \,b_n$, with $b_n \in \mathcal{B}_{13}$ for $i=2,13$ and
$b_n \in \mathcal{C}_{13}$ for $i=3$.
When this is the case, the list $L$ is processed to sum up the coefficients for the same $b_n$, that is,
if  $\lambda \,b_n, \lambda' \,b_n \in L$ for some $b_n \in \mathcal{C}_{13}$, then the two elements are removed from $L$ and
replaced by $(\lambda + \lambda')\,b_n$ only if $\lambda + \lambda' \neq 0$. For example, if in the end the list $L$ is $[(ab), (a^2c)st, (-a^2c)st, (bc), (bc^2)s]$, it will be processed for ``cleaning" the coefficients and become $L = [(ab + bc), (bc^2)s]$. 

We describe now what happens inside each iteration; this is the core of the algorithm.
 For every element $\alpha \in L$:
\begin{itemize}
\item[S1.] We check whether $\alpha \in \mathcal{C}_{13}$. If this is the case, there is nothing to do.	 \smallbreak
\item[S2.] If $\alpha \notin \mathcal{C}_{13}$, we check  whether $\alpha$ appears in one of the special cases. If this happens, we replace inside $L$ the element $\alpha$  by its summands given by this special case and we use the inverse Hecke relations to transform any negative powers appearing in the special case into positive ones before replacing $\alpha$. \smallbreak
\item[S3.] If $\alpha \notin \mathcal{C}_{13}$ and it does not appear in any special case, we should be able to apply some positive Hecke relation. We replace inside $L$ the element $\alpha$ with its summands arising from the positive Hecke relation. \smallbreak
\end{itemize}

The outputs of the \texttt{C++} program are the following: \smallbreak
\begin{itemize}
\item[O1.] The coefficients $\lambda_{j,n}^u\in \mathfrak{R}$ such that  $ub_{j}=\sum_{n=1}^{96}\lambda_{j,n}^{u}b_{n}$, for every $j\in\{1,\ldots, 96\}$.\smallbreak
	\item[O2.] The coefficients $\lambda_{j,n}^t\in \mathfrak{R}$ such that  $tb_{j}=\sum_{n=1}^{96}\lambda_{j,n}^{t}b_{n}$, for every $j\in\{1,\ldots, 96\}$.\smallbreak
		\item[O3.] The coefficients $\widetilde{\lambda}_{j,n}^s\in \mathfrak{R}$ such that  $sb_{j}=\sum_{n=1}^{100}\widetilde{\lambda}_{j,n}^{s}b_{n}$, for every
		$j\in\{1,\ldots,96\}\setminus\{81,93\}$.
		\smallbreak
\end{itemize}

\subsubsection{Second step: The products $sb_{81}$ and $sb_{93}$} From now on, we switch to the programming language SAGE \cite{sagemath} for all subsequent calculations (otherwise we would have needed to create programs in \texttt{C++} to do what SAGE can already do).
In order to proceed, 
we first need to express $sb_{81}$ and $sb_{93}$ as $\mathfrak{R}$-linear combinations of the elements of $\mathcal{C}_{13}$.
We have:
\smallbreak\smallbreak\smallbreak
\begin{center}
$\begin{array}[t]{lcl}
s\cdot b_{81}&=& s\cdot z^3tus\smallbreak\\
&=&z^3(stust)t^{-1}\smallbreak\\
&=&z^3(ustu)st^{-1}
\smallbreak\\
&=&z^3tustst^{-1}\smallbreak\\
&=&z^3tusts\,(d^{-1}t-cd^{-1})\smallbreak\\
&=&d^{-1}tuz^3stst-cd^{-1}tuz^3sts\smallbreak\\
&=&d^{-1}t\left(ub_{94}\right)-cd^{-1}t\left(ub_{82}\right)
\smallbreak\smallbreak\\
&=&d^{-1}t\sum_{\ell=1}^{96}\lambda_{94,\ell}^ub_\ell-d^{-1}ct\sum_{\ell=1}^{96}\lambda_{82,\ell}^ub_\ell\smallbreak\smallbreak\\
&=&d^{-1}\sum_{\ell=1}^{96}(\lambda_{94,\ell}^u-c\lambda_{82,\ell}^u)\,tb_\ell\smallbreak\smallbreak\\
&=&d^{-1}\sum_{n=1}^{96}\sum_{\ell=1}^{96}(\lambda_{94,\ell}^u-c\lambda_{82,\ell}^u)\,\lambda_{\ell,n}^t b_n\smallbreak\smallbreak\\
\end{array}$
$\begin{array}[t]{lcl}
s\cdot b_{93}&=& s\cdot z^3tust\smallbreak\\&=&z^3(stust)\smallbreak\\&=&z^3(ustu)s\smallbreak\\&=&z^3tusts\smallbreak\\&=&tuz^3sts
\smallbreak\\&=&t(u\cdot b_{82})
\smallbreak\\&=&t\sum_{\ell=1}^{96}\lambda_{82,\ell}^ub_\ell\smallbreak\\&=&\sum_{\ell=1}^{96}\lambda_{82,\ell}^utb_\ell\smallbreak\smallbreak\\&=&
\sum_{n=1}^{96}\sum_{\ell=1}^{96}\lambda_{82,\ell}^u\lambda_{\ell,n}^tb_n
\end{array}$
\end{center}
Therefore, if we set \smallbreak
\begin{center}
$
\widetilde{\lambda}_{81,n}^{s}:=d^{-1}\sum_{\ell=1}^{96}
(\lambda_{94,\ell}^u -c \lambda_{82,\ell}^u)\lambda_{\ell,n}^t\;\;\text{ and } \;\;
\widetilde{\lambda}_{93,n}^{s}:=\sum_{\ell=1}^{96}\lambda_{82,\ell}^u\lambda_{\ell,n}^t\,,$\quad for $n \in \{1,\ldots,96\}$,
\end{center}
and take $\widetilde{\lambda}_{81,n}^{s}=\widetilde{\lambda}_{93,n}^{s}=0$ for $n \in \{97,98,99,100\}$,
then we have 
\begin{center}
$sb_{j}=\sum_{n=1}^{100}\widetilde{\lambda}_{j,n}^{s}b_{n}$,\quad for every
		$j\in\{1,\ldots,96\}$.
\end{center}		

\subsubsection{Third step: Removing the elements $b_{97}$,  $b_{98}$,  $b_{99}$ and
 $b_{100}$}\label{tolabel}
Now, for every $j\in\{1,\ldots,96\}$, we will write $sb_j$ as an $\mathfrak{R}$-linear combination of elements in $\mathcal{B}_{13}$ instead of elements in $\mathcal{C}_{13}$. We are looking thus for the coefficients $\lambda_{j,n}^s$ such that
$sb_j =\sum_{n=1}^{96} \lambda_{j,n}^sb_n$.
For this, we need first to express the elements $b_{97}$, $b_{98}$, $b_{99}$, and $b_{100}$ as linear combinations of elements of  $\mathcal{B}_{13}$. 
We have:
$$\begin{array}{l}
b_{97}=tstsu=t(sb_8)=t\sum_{\ell=1}^{96}\lambda_{8,\ell}^sb_\ell=\sum_{\ell=1}^{96}\lambda_{8,\ell}^stb_\ell=
\sum_{n=1}^{96}\sum_{\ell=1}^{96}\lambda_{8,\ell}^s\lambda_{\ell,n}^tb_n
\\ \\
b_{98}=tstsut
=t(sb_{20})
=t\sum_{\ell=1}^{96}\lambda_{20,\ell}^sb_\ell
=\sum_{\ell=1}^{96}\lambda_{20,\ell}^stb_\ell
=\sum_{n=1}^{96}\sum_{\ell=1}^{96}\lambda_{20,\ell}^s\lambda_{\ell,n}^tb_n
\\ \\
b_{99}=utzuts
=u(tb_{36})
=u\sum_{\ell=1}^{96}\lambda_{36,\ell}^tb_\ell
=\sum_{\ell=1}^{96}\lambda_{36,\ell}^tub_\ell
=\sum_{n=1}^{96}\sum_{\ell=1}^{96}\lambda_{36,\ell}^t\lambda_{\ell,n}^ub_n
\\ \\
b_{100}=utzutst
=u(tb_{48})
=u\sum_{\ell=1}^{96}\lambda_{48,\ell}^tb_\ell
=\sum_{\ell=1}^{96}\lambda_{48,\ell}^tub_\ell
=\sum_{n=1}^{96}\sum_{\ell=1}^{96}\lambda_{48,\ell}^t\lambda_{\ell,n}^ub_n
\end{array}$$
Therefore, if we set
\begin{center}
$d_{97}^n:=\sum_{\ell=1}^{96}\lambda_{8,\ell}^s\lambda_{\ell,n}^t$,\;\; 
$d_{98}^n:=\sum_{\ell=1}^{96}\lambda_{20,\ell}^s\lambda_{\ell,n}^t$, 
\;\;$d_{99}^n:=\sum_{\ell=1}^{96}\lambda_{36,\ell}^t\lambda_{\ell,n}^u$,\;\; and \;\;$d_{100}^n:=\sum_{\ell=1}^{96}\lambda_{48,\ell}^t\lambda_{\ell,n}^u$,
\end{center}
then we have
\begin{center}
$b_i = \sum_{n=1}^{96} d_i^n b_n$, \quad for all $i=97,98,99,100$.
\end{center}

We notice now the presence of $\lambda_{8,\ell}^s$ and $\lambda_{20,\ell}^s$, $\ell=1,\ldots,96$, in  $d_{97}^n$ and $d_{98}^n$. We will calculate these coefficients by examining the cases of $sb_8$ and $sb_{20}$ separately.
According to the \texttt{C++} program, we have 
$$\widetilde{\lambda}_{8,97}^{s}=\widetilde{\lambda}_{8,98}^{s}=
\widetilde{\lambda}_{8,100}^{s}=0
\quad\text{and}\quad
\widetilde{\lambda}_{20,97}^{s}=\widetilde{\lambda}_{20,98}^{s}=
\widetilde{\lambda}_{20,99}^{s}=0.$$
Therefore, we obtain
\begin{center}
$sb_{8}=\sum_{n=1}^{96}\widetilde{\lambda}_{8,n}^{s}b_{n}+
\widetilde{\lambda}_{8,99}^{s}\sum_{n=1}^{96}d^n_{99}b_n\quad\text{and}\quad sb_{20}=\sum_{n=1}^{96}\widetilde{\lambda}_{20,n}^{s}b_{n}+
\widetilde{\lambda}_{20,100}^{s}\sum_{n=1}^{96}d_{100}^nb_n,$
\end{center}
whence 
\begin{equation}\label{l8l20}
\lambda_{8,n}^s=\widetilde{\lambda}_{8,n}^{s}+\widetilde{\lambda}_{8,99}^{s}d^n_{99}
\quad\text{and}\quad 
\lambda_{20,n}^s=\widetilde{\lambda}_{20,n}^{s}+\widetilde{\lambda}_{20,100}^{s}d^n_{100}, \quad \text{for all } n \in \{1,\ldots,96\}.
\end{equation}

We can now calculate all remaining $\lambda_{j,n}^s$.
For every $j \in \{1,\ldots,96\} \setminus \{8,20\}$, we have
$$\begin{array}[t]{lcl}
sb_j
&=&\sum_{n=1}^{96}\widetilde{\lambda}_{j,n}^{s}b_{n}
+\sum_{i=97}^{100}\widetilde{\lambda}_{j,i}^{s}\sum_{n=1}^{96}d_{i}^nb_n
\smallbreak\smallbreak\\
\end{array}$$
Therefore, for every $j \in \{1,\ldots,96\} \setminus \{8,20\}$, we obtain
\begin{equation}\label{li}
\lambda_{j,n}^s=\widetilde{\lambda}_{j,n}^{s}+
\sum_{i=97}^{100} \widetilde{\lambda}_{j,i}^{s} d_{i}^n.
\end{equation}

To summarise, from the outputs of the \texttt{C++} program, and from Equations \eqref{l8l20} and \eqref{li}, we have the expressions of $sb_j$, $tb_j$ and $ub_j$ as $\mathfrak{R}$-linear combinations of elements of $\mathcal{B}_{13}$, for every $j\in\{1,\dots,96\}$.

\subsubsection{Fourth step: The inductive procedure} 
Let  $k\in\{0,1,2,3\}$. We observe that: 
\begin{equation}\label{ind}
  \begin{array}{lll}
  b_{24k+2\phantom{1}}= b_{24k+1\phantom{3}}\cdot u,\;&
  b_{24k+3\phantom{1}}=  b_{24k+1\phantom{1}}\cdot s,&
  \;b_{24k+13}=  b_{24k+1\phantom{1}}\cdot t,\smallbreak\smallbreak\\
  b_{24k+4\phantom{1}}= b_{24k+13}\cdot s,&
  b_{24k+5\phantom{1}}= b_{24k+3\phantom{1}}\cdot u,&
  \;b_{24k+6\phantom{1}}= b_{24k+2\phantom{1}}\cdot s,\smallbreak\smallbreak\\
  b_{24k+7\phantom{1}}= b_{24k+13}\cdot u,&
  b_{24k+8\phantom{1}}= b_{24k+4\phantom{1}}\cdot u,&
  \;b_{24k+9\phantom{1}}= b_{24k+7\phantom{1}}\cdot s,\smallbreak\smallbreak\\
  b_{24k+15}= b_{24k+3\phantom{1}}\cdot t,&
  b_{24k+10}= b_{24k+15}\cdot s,&
  \;b_{24k+11}= b_{24k+15}\cdot u,\smallbreak\smallbreak\\
  b_{24k+14}= b_{24k+2\phantom{1}}\cdot t,&
  b_{24k+12}= b_{24k+14}\cdot s,&\;
  b_{24k+m}=b_{24k+m-12}\cdot t, \;\;m\in\{16,\ldots, 24\}.
  \end{array}
  \end{equation}

Thus, for every $m \in \{2,\ldots,24\}$, there exists a unique $m' \in \{1,\ldots,15\}$
and $g \in \{s,t,u\}$ such that
$$ b_{24k+m} = b_{24k+m'} \cdot g.$$
Then, for all $j \in \{1,\ldots,96\}$, we have
$$ b_{24k+m}b_j =   b_{24k+m'} \sum_{n=1}^{96} \lambda_{j,n}^gb_n, $$
whence
\begin{equation}\label{eq3}
 \tau(b_{24k+m}b_j) =   \sum_{n=1}^{96} \lambda_{j,n}^g\tau(b_{24k+m'} b_n).
 \end{equation}
As a consequence, we can compute the entries of the matrix $A =(\tau(b_ib_j))_{1 \leq i,j \leq 96}$
as follows: We start with $k=0$.
We first complete the row $24k+1$. 
By definition of $\tau$, we have $\tau(b_1b_j) = \tau(b_j)=\delta_{1j}$, for all $j \in \{1,\ldots,96\}$. 
Then, thanks to \eqref{eq3}, 
we can complete the entries of the row $24k+m$, for $m \in \{2,\ldots,24\}$, if we have already completed the entries of the row $24k+m'$, with $m'$ and $g$ given by  \eqref{ind}. 
This is why, we complete the entries of the rows of $A$ in the following order (we just give the number of the row): $24k+2$,  $24k+3$, $24k+13$,
$24k+4$, $24k+5$, $24k+6$, $24k+7$, $24k+8$, $24k+9$, 
$24k+15$, $24k+10$, $24k+11$, $24k+14$, $24k+12$,
$24k+16$, $24k+17$, $24k+18$, $24k+19$,
$24k+20$, $24k+21$, $24k+22$, $24k+23$,   $24k+24$.
We repeat the same procedure for $k =1,2,3$, but first we have to show how we fill in the entries of the row $24k+1$ for $k \neq 0$.

Assume $k \in \{1,2,3\}$. We have $b_{24k+1}=z^k$. 
We distinguish two cases:
\begin{itemize}
	\item If $1 \leq j \leq 24(4-k)$, then $b_{24k+1}b_j=z^kb_j$.  Since $24k+1\leq 24k+j\leq 96$, we have $z^kb_j=b_{24k+j}$. We have assumed that $k\not=0$, so $b_{24k+j}\in \mathcal{B}_{13} \setminus \{1\}$.
	Hence, 
	\begin{equation}\label{eq2}
	\tau(b_{24k+1}b_j )=0.
	\end{equation}
	\item If $24(4-k) < j \leq 96$, then $b_{24k+1}b_j=z^{k}b_j=
	z^{k-4} b_{j}z^4 = b_{24k+j-96}\cdot z^4$. If $	z^4=\sum_{n=1}^{96}\mu_{n}b_{n}$ with $\mu_n \in \mathfrak{R}$, then
	\begin{equation}\label{eq4}
	\tau(b_{24k+1}b_j)=\tau(b_{24k+j-96}\cdot z^4)=\sum_{n=1}^{96}\mu_{n}
	\tau(b_{24k+j-96}\,b_{n}).
	\end{equation} 
\end{itemize}

Notice now that the values $\tau(b_{24k+j-96}\,b_{n})$ in \eqref{eq4} are entries of row
	$24k+j-96$ of the matrix $A$. Since $24k+j-96 < 24k+1$, these are entries that have already been computed when we try to complete row $24k+1$.

So, for the algorithm to be complete, it remains to calculate the coefficients $\mu_{n}$ of $z^4$ when expressed as an $\mathfrak{R}$-linear combination of elements of $\mathcal{B}_{13}$. We have:
$$\begin{array}[t]{lcl}
z^4&=&z^3ustustust\smallbreak\smallbreak\\
&=&ustu(s\cdot z^3tust)\smallbreak\smallbreak\\
&=&ustu\sum_{\ell}\lambda_{93,\ell}^sb_\ell\smallbreak\smallbreak\\
&=&ust\sum_{\ell,p}\lambda_{93,\ell}^s\lambda_{\ell,p}^ub_p\smallbreak\smallbreak\\
&=&us\sum_{\ell,p,q}\lambda_{93,\ell}^s\lambda_{\ell,p}^u\lambda_{p,q}^tb_q\smallbreak\smallbreak\\
&=&u\sum_{\ell,p,q,r}\lambda_{93,\ell}^s\lambda_{\ell,p}^u\lambda_{p,q}^t
\lambda_{q,r}^sb_r
\smallbreak\smallbreak\\
&=&\sum_{\ell,p,q,r,n}\lambda_{93,\ell}^s\lambda_{\ell,p}^u\lambda_{p,q}^t
\lambda_{q,r}^s\lambda_{r,n}^ub_n
\smallbreak\smallbreak\\
\end{array}$$
Therefore, for every $n\in\{1,\ldots,96\}$, we have
$\mu_{n}=\sum_{\ell,p,q,r}\lambda_{93,\ell}^s\lambda_{\ell,p}^u\lambda_{p,q}^t
\lambda_{q,r}^s\lambda_{r,n}^u.$

\subsection{Proof of the third condition}
We will prove here Condition (3) of Conjecture \ref{BMM sym}. 
Given that our basis $\mathcal{B}_{13}$ satisfies the properties required by the lifting conjecture of Malle and Michel, 
 it is enough to prove that $\tau$ satisfies Equation \eqref{extra con2}, that is, we have to show that 
  $\tau(T_{\boldsymbol{\pi}}b_i^{-1})=0$, for all $i=2,\ldots,96$. Since   $T_{\boldsymbol{\pi}}=T_{\boldsymbol{\beta}}^{|Z(G_{13})|}=z^4$, we will show that $\tau(z^4b_i^{-1})=0$, for all $i=2,\ldots,96$.
  We will make use of the fact that $\tau$ is a trace function, since we proved it in the previous subsection.
 
 We first consider the case $i>24$. We write $b_i$ as $z^kb_m$, with $k \in \{1,2,3\}$ and $m \in \{1,2,\ldots,24\}$.  We have $z^4b_i^{-1}=z^{4-k}b_m^{-1}$. If $m\in\{13,\ldots,24\}$, then $b_m=b_{m-12}t$, and so  
 $$\tau(z^{4-k}b_m^{-1})=\tau(z^{4-k}t^{-1}b_{m-12}^{-1})=\tau(z^{4-k}b_{m-12}^{-1}t^{-1})=d^{-1}\tau(z^{4-k}b_{m-12}^{-1}t)-cd^{-1}\tau(z^{4-k}b_{m-12}^{-1}).$$ 
As a result,  we can restrict ourselves to proving that  $\tau(z^{4-k}b_{m}^{-1})=\tau(z^{4-k}b_{m}^{-1}t)=0$ for every $m\in\{1,\ldots,12\}$.
Using the inverse Hecke relations we can write every $b_m^{-1}$, $m=1,\ldots,12$, as a linear combination of $b_j$, $j=1,\ldots,24$. Since $\tau$ is linear, it suffices to show that $\tau(z^{4-k}b_j)=\tau(z^{4-k}b_jt)=0$, for every $j\in\{1,\ldots,24\}$.
Since $4-k\in\{1,2,3\}$, we have that $z^{4-k}b_j\in\{b_{25},\ldots,b_{96}\}\subset \mathcal{B}_{13}\setminus\{1\}$ and, hence, $\tau(z^{4-k}b_j)=0$.
In order to prove that $\tau(z^{4-k}b_jt)=0$, we distinguish the following two cases:
For every $j\in\{1,\ldots,12\}$ we have $z^{4-k}b_jt=z^{4-k}b_{12+j}\in\mathcal{B}_{13}\setminus\{1\}$ and, hence, $\tau(z^{4-k}b_jt)=0$.
	 For every $j\in\{13,\ldots,24\}$ we have $z^{4-k}b_jt=z^{4-k}b_{j-12}t^2=cz^{4-k}b_{j-12}t+dz^{4-k}b_{j-12}=cz^{4-k}b_j+dz^{4-k}b_{j-12}$. As a result, $z^{4-k}b_jt$ is a linear combination of elements of $\{b_{25},\ldots,b_{96}\}\subset \mathcal{B}_{13}\setminus\{1\}$ and, hence, $\tau(z^{4-k}b_jt)=0$.

Now, let us consider the case $i \leq 24$. We recall that $z=(stu)^3=(tus)^3=(ust)^3$. Using the fact that $\tau$ is a trace function and that the entries of the matrix $A$ are known, we have{\small
	$$\begin{array}{rcl}
	\tau(z^4b_2^{-1}) &=& d^{-1}\tau(z^4u)-cd^{-1}\tau(z^4) = d^{-1}\left(\tau(b_{73}b_{26})
	-c\tau(b_{73}b_{25})\right)=d^{-1}(cb^6d^{12} - cb^6d^{12})=0\smallbreak\smallbreak\\
	\tau(z^4b_3^{-1}) &=& b^{-1}\tau(z^4s)-ab^{-1}\tau(z^4) = b^{-1}\left(\tau(b_{73}b_{27})
	-a\tau(b_{73}b_{25})\right)=b^{-1}(ab^6d^{12} - ab^6d^{12})=0\smallbreak\smallbreak\\
	\tau(z^4b_4^{-1}) &=& \tau(z^3tustustuss^{-1}t^{-1})=\tau(z^3tust(ustu)t^{-1})=\tau(z^3tust\cdot tus)=\tau(b_{93}b_{9})=0\smallbreak\smallbreak\\
	\tau(z^4b_5^{-1}) &=& \tau(z^3stustustuu^{-1}s^{-1})= \tau(z^3stustusts^{-1})=\tau(z^3s^{-1}stustust)=\tau(z^3tust\cdot ust)=
	\tau(b_{93}b_{18})=0\smallbreak\smallbreak\\
	\tau(z^4b_6^{-1}) &=& \tau(z^4s^{-1}u^{-1})=\tau(z^4u^{-1}s^{-1})=\tau(z^4b_5^{-1})=0\smallbreak\smallbreak\\
	\tau(z^4b_7^{-1}) &=& \tau(z^3stustustuu^{-1}t^{-1})=\tau(z^3stustus)=
	\tau(z^3s^2tustu) =
	a\tau(z^3stu\cdot stu)+b \tau(z^3tu\cdot stu)=\smallbreak\smallbreak\\&=&a\tau(b_{83}b_{11})+b\tau(b_{79}b_{11})=0+0=0\smallbreak\smallbreak\\
	\tau(z^4b_8^{-1}) &=& \tau(z^3stu(stust)uu^{-1}s^{-1}t^{-1})=\tau(z^3stu(ustu)ss^{-1}t^{-1})
	=\tau(z^3stutustt^{-1})=\tau(z^3stu\cdot tus)=\smallbreak\smallbreak\\&=&\tau(b_{83}b_{9})=0\smallbreak\smallbreak\\
	\tau(z^4b_9^{-1}) &=& \tau(z^3tustustuss^{-1}u^{-1}t^{-1})=\tau(z^3tus\cdot tus)=\tau(b_{81}b_{9})=0\smallbreak\smallbreak\\
	\tau(z^4b_{10}^{-1}) &=& \tau(z^3tust(ustu)ss^{-1}t^{-1}s^{-1})
	=\tau(z^3tusttustt^{-1}s^{-1})=\tau(z^3tust\cdot tu)=\tau(b_{93}b_7)=0\smallbreak\smallbreak\\
	\tau(z^4b_{11}^{-1})&=&\tau(z^4u^{-1}t^{-1}s^{-1})=\tau(z^4s^{-1}u^{-1}t^{-1})
	=\tau(z^4b_9^{-1})=0\smallbreak\smallbreak\\
	\tau(z^4b_{12}^{-1})&=& \tau(z^4s^{-1}t^{-1}u^{-1})=\tau(z^4u^{-1}s^{-1}t^{-1})
	=\tau(z^4b_8^{-1})=0\smallbreak\smallbreak\\
	\tau(z^4b_{13}^{-1})&=& d^{-1}\tau(z^4t)-cd^{-1}\tau(z^4) = d^{-1}\left(\tau(b_{73}b_{37})-c\tau(b_{73}b_{25})\right)=d^{-1}(cb^6d^{12} - cb^6d^{12})=0\smallbreak\smallbreak\\
		\tau(z^4b_{14}^{-1})&=& \tau(z^4t^{-1}u^{-1})=\tau(z^4u^{-1}t^{-1})=\tau(z^4b_7^{-1})=0\smallbreak\smallbreak\\
	\tau(z^4b_{15}^{-1})&=& \tau(z^4t^{-1}s^{-1})=\tau(z^4s^{-1}t^{-1})=\tau(z^4b_4^{-1})=0\smallbreak\smallbreak\\
	\tau(z^4b_{16}^{-1})&=& \tau(z^3ust(ustu)stt^{-1}s^{-1}t^{-1})=
	\tau(z^3usttustt^{-1})=\tau(z^3ust\cdot tus)=\tau(b_{90}b_9)=0\smallbreak\smallbreak\\
	\tau(z^4b_{17}^{-1})&=& \tau(z^4t^{-1}u^{-1}s^{-1})=\tau(z^4s^{-1}t^{-1}u^{-1})=\tau(z^4b_{12}^{-1})=0\smallbreak\smallbreak\\
	\tau(z^4b_{18}^{-1})&=& \tau(z^4t^{-1}s^{-1}u^{-1})=\tau(z^4u^{-1}t^{-1}s^{-1})=\tau(z^4b_{11}^{-1})=0\smallbreak\smallbreak\\
	\tau(z^4b_{19}^{-1})&=& \tau(z^4t^{-1}u^{-1}t^{-1}) = \tau(z^4t^{-2}u^{-1})=
	d^{-1}\tau(z^4u^{-1})-cd^{-1} \tau(z^4t^{-1}u^{-1}) =\smallbreak\smallbreak\\
	&=&d^{-1}\tau(z^4b_2^{-1})-cd^{-1}\tau(z^4b_{14}^{-1})=0\smallbreak\smallbreak\\
	\tau(z^4b_{20}^{-1})&=& \tau(z^4t^{-1}u^{-1}s^{-1}t^{-1}) = \tau(z^4t^{-2}u^{-1}s^{-1})=
	d^{-1}\tau(z^4u^{-1}s^{-1})-cd^{-1} \tau(z^4t^{-1}u^{-1}s^{-1}) =\smallbreak\smallbreak\\
	&=&
	d^{-1}\tau(z^4b_5^{-1})-cd^{-1} \tau(z^4b_{17}^{-1})=0\smallbreak\smallbreak\\
	\tau(z^4b_{21}^{-1})&=& \tau(z^4t^{-1}s^{-1}u^{-1}t^{-1}) = \tau(z^4t^{-2}s^{-1}u^{-1})= 
	d^{-1}\tau(z^4s^{-1}u^{-1})-cd^{-1} \tau(z^4t^{-1}s^{-1}u^{-1}) =\smallbreak\smallbreak\\
	&=&d^{-1}\tau(z^4b_6^{-1})-cd^{-1} \tau(z^4b_{18}^{-1})=0\smallbreak\smallbreak\\
	\tau(z^4b_{22}^{-1})&=&\tau(z^3ust(ustu)stt^{-1}s^{-1}t^{-1}s^{-1})=
	\tau(z^3usttustt^{-1}s^{-1})=\tau(z^3ust\cdot tu)=\tau(b_{90}b_7)=0\smallbreak\smallbreak\\
	\tau(z^4b_{23}^{-1})&=&\tau(z^3(ustus)tustt^{-1}u^{-1}t^{-1}s^{-1})=
	\tau(z^3u^{-1}t^{-1}s^{-1}stusttus)=\tau(z^3st\cdot tus)=\tau(b_{87}b_9)=0\smallbreak\smallbreak\\
	\tau(z^4b_{24}^{-1})&=&\tau(z^4t^{-1}s^{-1}t^{-1}u^{-1})=
	\tau(z^4t^{-1}u^{-1}t^{-1}s^{-1})=\tau(z^4b_{23}^{-1})=0.
	\end{array}$$}

\subsection{Further results and remarks}	

We have now shown that the linear map $\tau : \mathcal{H}(G_{13}) \rightarrow R(G_{13})$ defined by $\tau(b)=\delta_{1b}$, for all $b \in \mathcal{B}_{13}$, is the canonical symmetrising trace on $\mathcal{H}(G_{13})$ and, thus, the BMM symmetrising trace conjecture is satisfied for all exceptional $2$-reflection groups of rank $2$. 
In the process, we have also proved the following:
\begin{thm}
The lifting conjecture  holds for $G_{13}$.
\end{thm}

In our proof, we have used the fact that the set $\mathcal{B}_{13}$ is a basis of   $\mathcal{H}(G_{13})$. We will now see that our proof
provides a proof of this fact. 

As remarked at the end of \S \ref{tolabel}, from the outputs of the \texttt{C++} program, and from Equations \eqref{l8l20} and \eqref{li}, we have the expressions of $sb_j$, $tb_j$ and $ub_j$ as $R(G_{13})$-linear combinations of elements of $\mathcal{B}_{13}$, for every $j\in\{1,\dots,96\}$. 
This in fact allows us to express any product of the generators $s$, $t$ and $u$, and in turn any element of $\mathcal{H}(G_{13})$, as a linear combination of the elements of $\mathcal{B}_{13}$. 
We give the following example:
$$\begin{array}{lcl}
sust&=&su\sum_{\ell }\lambda^s_{13,\ell}b_\ell = s\sum_{\ell,r}\lambda^s_{13,l}
\lambda^u_{\ell,r}b_r=\sum_{\ell,r,n}\lambda^s_{13,l}\lambda^u_{\ell,r}\lambda^s_{r,n}b_n.
\end{array}$$

\begin{rem}\rm  Let $g \in \{s,t,u\}$ and $m \in \Z^*$. An interesting observation about $G_{13}$ is that, due to the fact that $G_{13}$ is a $2$-reflection group, and thanks to the inverse Hecke relations, $g^m$ can be always written as a linear combination of $g$ and $1$. This in turn implies that any product of positive or negative powers of the generators can be expressed as a linear combination of products where these powers are either $0$ or $1$.
\end{rem}

Therefore, we can deduce the following:

\begin{prop}
The set $\mathcal{B}_{13}$ is a spanning set for $\mathcal{H}(G_{13})$ as an $R(G_{13})$-module.
\end{prop}

Since $|\mathcal{B}_{13}| = |G_{13}|=96$, Theorem \ref{mention} implies the following:

\begin{thm}\label{Bn basis}
 The set $\mathcal{B}_{13}$ is a basis for  $\mathcal{H}(G_{13})$ as an $R(G_{13})$-module. In particular, the BMR freeness conjecture holds for $G_{13}$.
\end{thm}

\begin{rem}\rm As we also remarked in \cite{BCCK}, this does not mean that we have obtained a computerised proof of the BMR freeness conjecture for $G_{13}$, because the ``special cases'' incorporated in the \texttt{C++} program are the exact calculations made  by hand by the second author for the proof of the BMR freeness conjecture in \cite{Ch17}.
\end{rem}

\section{Appendix}

We will discuss here the 25 special cases mentioned in \S\ref{first step}, and give examples of their use. Let $g \in \{s,t,u\}$ and $b_j \in \mathcal{B}_{13}$.

Cases 1--3 establish that $z$ is central. More precisely, for every $k\in\{1,2,3\}$ we have $g\cdot z^k \cdot b_j=z^k \cdot g \cdot b_j$. 
We created these cases, because otherwise the computer took a lot of time, applying the braid relations over and over, simply to commute $z$.

Cases 4--23  express a given $gb_j$  as a sum of other elements in $\mathcal{H}(G_{13})$, which
	 can be in turn expressed as  linear combinations of elements in $\mathcal{B}_{13}$ if 
	 we apply  the inverse Hecke relations, the positive Hecke relations, or/and some  other special case(s). In order to make it clearer to the reader, we give the following example, where we express $tb_{58}=tz^2sts$ as an $\mathfrak{R}$-linear combination of elements in $\mathcal{B}_{13}$. We first use the fact that $z$ is central and we write $tz^2sts=z^2tsts$. 
We then apply
$ $\\ \\
{\bf Case 8: } For $k\in\{2,3\}$ and $l\in\{0,1\}$, we have
$$\begin{array}{lcl}
z^ktstst^l &=& {c}z^{k}stst^l+{ad}z^{k}st^l+{bcd}z^{k}t^{l-1}+
{abd^2}z^{k}t^{-1}s^{-1}t^{l-1}+{b^2d^2}z^{k-1}ust^2ut^{l}.\end{array}$$
$ $\\
We obtain:
	$$\begin{array}{lcl}z^2tsts
	&=&c\boldsymbol{b_{58}}+ad\boldsymbol{b_{51}}+bcdz^2t^{-1}+abd^2z^2t^{-1}s^{-1}t^{-1}+b^2d^2zust^2u.\smallbreak
	\end{array}$$
		We now apply  the inverse Hecke relations and we get:
		$$\begin{array}{ccl}
	bcdz^2t^{-1}&=&bcz^2t-c^2dz^2=bc\boldsymbol{b_{61}}-c^2b\boldsymbol{b_{49}}\smallbreak\smallbreak\smallbreak\smallbreak\\
	abd^2z^2t^{-1}s^{-1}t^{-1}&=& -a^2c^2z^2+a^2cz^2t+ac^2z^2s-acz^2ts+a^2cz^2t-acz^2
	st+az^2tst-a^2z^2t^2\smallbreak\smallbreak\\
	&=&-a^2c^2\boldsymbol{b_{49}}+a^2c\boldsymbol{b_{61}}+ac^2\boldsymbol{b_{51}}
	-ac\boldsymbol{b_{52}}+a^2c\boldsymbol{b_{61}}-ac\boldsymbol{b_{63}}
	+a\boldsymbol{b_{64}}-a^2z^2t^2\smallbreak\smallbreak\\
	\end{array}$$
By the positive Hecke relation for $t$,  we have:
\smallbreak\smallbreak
	$\begin{array}{rcl}-a^2z^2t^2\;=\;-a^2cz^2t-a^2dz^2\;=\;-a^2c\boldsymbol{b_{61}}-a^2d\boldsymbol{b_{49}}.\smallbreak\smallbreak\smallbreak
	\end{array}$
	\newline
	Hence,
	\smallbreak\smallbreak
	$\begin{array}{lcl}
	abd^2z^2t^{-1}s^{-1}t^{-1}&=&-a^2c^2\boldsymbol{b_{49}}+a^2c\boldsymbol{b_{61}}+ac^2\boldsymbol{b_{51}}
	-ac\boldsymbol{b_{52}}+a^2c\boldsymbol{b_{61}}-ac\boldsymbol{b_{63}}
	+a\boldsymbol{b_{64}}-\smallbreak\\&&-a^2c\boldsymbol{b_{61}}-a^2d\boldsymbol{b_{49}}
	\smallbreak\smallbreak\\
	&=&(-a^2c^2-a^2d)\boldsymbol{b_{49}}+ac^2\boldsymbol{b_{51}}
	-ac\boldsymbol{b_{52}}+a^2c\boldsymbol{b_{61}}-ac\boldsymbol{b_{63}}
	+a\boldsymbol{b_{64}}
	\smallbreak\smallbreak
	\end{array}$
	\newline
	
	It remains to write the element $b^2d^2zust^2u$ as an $\mathfrak{R}$-linear combination of elements in $\mathcal{B}_{13}$.
	We apply again the positive Hecke relation for $t$ and we have
	$$b^2d^2zust^2u=b^2d^2czustu+b^2d^3zusu$$
	We now apply
$ $\\ \\	
{\bf Case 12: } For $k\in\{0,1,2,3\}$ and $l\in\{0,1\}$, we have $z^kustut^l=z^ktust^{l+1}$
$ $\\  \\ and \\. \\
{\bf Case 14:} For $k\in\{1,2,3\}$ and $l\in\{0,1\}$, we have
$$\begin{array}{lcl} z^kusut^l&=&{c}z^{k}ust^l+
{ad}z^{k}t^l+{bcd}z^{k}s^{-1}u^{-1}t^l+
{bd^2}z^{k-1}stu^2st^{l+1}\end{array}
$$
to obtain:	
	 $$\begin{array}{ccl}
	b^2d^2czustu&=&b^2d^2cztust=b^2d^2c\boldsymbol{b_{45}}\smallbreak\smallbreak\smallbreak\smallbreak\\
	b^2d^3zusu&=&b^2d^3czus+ab^2d^4z+b^3cd^4zs^{-1}u^{-1}+b^3d^5stu^2st\smallbreak\\
	&=&b^2d^3c\boldsymbol{b_{30}}+ab^2d^4\boldsymbol{b_{25}}+b^3cd^4zs^{-1}u^{-1}+b^3d^5stu^2s
	\end{array}$$ 
	We use the inverse Hecke relations for $s$ and $u$ and the positive Hecke relation for $u$ and we get:
	\smallbreak\smallbreak
	$\begin{array}{ccl}
	b^3cd^4zs^{-1}u^{-1}&=&b^2d^3ac^2z-acb^2d^3zu-c^2b^2d^3zs+czsu\smallbreak\\
	&=&	b^2d^3ac^2\boldsymbol{b_{25}}-acb^2d^3\boldsymbol{b_{26}}-c^2b^2d^3\boldsymbol{b_{27}}+cb^2d^3\boldsymbol{b_{29}}\smallbreak
	\smallbreak\smallbreak\\
	b^3d^5stu^2st&=&b^3cd^5stust+b^3d^6stst\smallbreak\\
	&=&b^3cd^5stust+b^3d^6\boldsymbol{b_{22}}\smallbreak\smallbreak\smallbreak
	\end{array}$
	\newline
	We apply
	$ $\\ \\
	{\bf Case 5:} Fore $k\in\{0,1, 2\}$ and $l\in\{0,1\}$, we have
$z^kstust^l=z^{k+1}t^{-1}s^{-1}u^{-1}t^{l-2}$
$ $\\ \\
and we obtain:
$$b^3cd^5stust=b^3cd^5zt^{-1}s^{-1}u^{-1}t^{-1}.$$
Using the inverse Hecke relations, and then the positive Hecke relations, 
we have:
		$$\begin{array}{lcl}
	b^3cd^5  stust&=&
	(d^3ab^2c^2 + c^4d^2ab^2)\boldsymbol{b_{25}}
 - c^3d^2ab^2\boldsymbol{b_{26}}
- c^4d^2b^2\boldsymbol{b_{27}}+
	d^2b^2c^3\boldsymbol{b_{28}}+
 c^3d^3b^2\boldsymbol{b_{29}}+\smallbreak\smallbreak\\
&&+d^2ab^2c^2\boldsymbol{b_{31}}
 - d^2b^2c^2\boldsymbol{b_{32}}
 - d^2ab^2c^3\boldsymbol{b_{37}}+
  c^2d^2ab^2\boldsymbol{b_{38}}
 + c^3d^2b^2\boldsymbol{b_{39}}
   - d^2b^2c^2\boldsymbol{b_{40}}-\smallbreak\smallbreak\\
	  && - c^2d^2b^2\boldsymbol{b_{41}}
	  - d^2ab^2c\boldsymbol{b_{43}}
	 + d^2b^2c\boldsymbol{b_{44}}
		\end{array}$$
		\smallbreak\smallbreak
	We conclude the following:
	\smallbreak\smallbreak
	$$\begin{array}{lcl}
	z^2tsts&=&c\boldsymbol{b_{58}}+ad\boldsymbol{b_{51}}+
	bc\boldsymbol{b_{61}}-bc^2\boldsymbol{b_{49}}+
	(-a^2c^2-a^2d)\boldsymbol{b_{49}}+ac^2\boldsymbol{b_{51}}
	-ac\boldsymbol{b_{52}}+a^2c\boldsymbol{b_{61}}+
	\smallbreak\smallbreak\\
	&&+
	ac\boldsymbol{b_{63}}
	+a\boldsymbol{b_{64}}
	
	+b^2d^2c\boldsymbol{b_{45}}+b^2d^3c\boldsymbol{b_{30}}+ab^2d^4\boldsymbol{b_{25}}+
	b^2d^3ac^2\boldsymbol{b_{25}}-acb^2d^3\boldsymbol{b_{26}}-\smallbreak\smallbreak\\
	&&-c^2b^2d^3\boldsymbol{b_{27}}+cb^2d^3\boldsymbol{b_{29}}+

	(d^3ab^2c^2 + c^4d^2ab^2)\boldsymbol{b_{25}}
	- c^3d^2ab^2\boldsymbol{b_{26}}
	- c^4d^2b^2\boldsymbol{b_{27}}+\smallbreak\smallbreak\\
	&&+d^2b^2c^3\boldsymbol{b_{28}}+
	c^3d^3b^2\boldsymbol{b_{29}}+
	d^2ab^2c^2\boldsymbol{b_{31}}
	- d^2b^2c^2\boldsymbol{b_{32}}
	- d^2ab^2c^3\boldsymbol{b_{37}}+
	c^2d^2ab^2\boldsymbol{b_{38}}+\smallbreak\smallbreak\\
	&&
	+ c^3d^2b^2\boldsymbol{b_{39}}
	- d^2b^2c^2\boldsymbol{b_{40}}
 - c^2d^2b^2\boldsymbol{b_{41}}
	- d^2ab^2c\boldsymbol{b_{43}}
	+ d^2b^2c\boldsymbol{b_{44}}
	
	+b^3d^6\boldsymbol{b_{22}}
	\smallbreak\smallbreak\smallbreak\smallbreak\\
	
	&=&b^3d^6\boldsymbol{b_{22}}+(d^3ab^2c^2 + c^4d^2ab^2+ab^2d^4+2b^2d^3ac^2)\boldsymbol{b_{25}}+
	(-acb^2d^3-c^3d^2ab^2)\boldsymbol{b_{26}}+\smallbreak\smallbreak\\
	&&(- d^3b^2c^2 - d^2b^2c^4)\boldsymbol{b_{27}}+
d^2b^2c^3\boldsymbol{b_{28}}+
	(c^3d^3b^2+cb^2d^3)\boldsymbol{b_{29}}+
	d^3b^2c\boldsymbol{b_{30}}+
	d^2ab^2c^2\boldsymbol{b_{31}}-\smallbreak\smallbreak\\&&
	- d^2b^2c^2\boldsymbol{b_{32}}
	- d^2ab^2c^3\boldsymbol{b_{37}}+
	c^2d^2ab^2\boldsymbol{b_{38}}
	+ c^3d^2b^2\boldsymbol{b_{39}}
	- d^2b^2c^2\boldsymbol{b_{40}}
	- c^2d^2b^2\boldsymbol{b_{41}}-\smallbreak\smallbreak\\
	&&
	- d^2ab^2c\boldsymbol{b_{43}}
	+ d^2b^2c\boldsymbol{b_{44}}+
	d^2b^2c\boldsymbol{b_{45}}+
	(-bc^2-a^2d-a^2c^2)\boldsymbol{b_{49}}+
	(da+ac^2)\boldsymbol{b_{51}}-\smallbreak\\
	&&-ac\boldsymbol{b_{52}}+
	c\boldsymbol{b_{58}}+(bc+a^2c)\boldsymbol{b_{61}}-ac\boldsymbol{b_{63}}+a\boldsymbol{b_{64}}\smallbreak\smallbreak
	\end{array}$$

Finally, Cases 24--25 express a given $sb_j$  as a sum of other elements in $\mathcal{H}(G_{13})$. As in the previous example, we use the Hecke relations and other special cases to write this time $sb_j$ as a linear combination of elements in  $\mathcal{C}_{13}:=\mathcal{B}_{13}\cup\{b_{97}:=tstsu,\;b_{98}:=tstsut,\;b_{99}:=zututs,\;b_{100}:=zututst\}$. The reason we introduce $\mathcal{C}_{13}$ is to avoid the creation of (many) more special cases. 
The four extra elements can be expressed as  $\mathfrak{R}$-linear combinations of elements in the basis, as detailed in \S \ref{tolabel}. They are four, and not two, because each special case is ``double'' in the sense that it involves a basis element ending in $t^l$ with $l \in \{0,1\}$. 
This comes from the fact that $\mathcal{E}_{13}^{k,1}=\mathcal{E}_{13}^{k,0}t$, for every $k\in\{0,1,2,3\}$.

 As we can see from the examples, doing all the calculations by hand is possible (it has already been done in \cite{Ch17}), but time-consuming. Moreover, simple mistakes can be made, like forgetting a coefficient. This is why we created the
  \texttt{C++} program described in \S\ref{first step}.


\begin{thebibliography}{99}
%
\addcontentsline{toc}{chapter}{Bibliography}
\bibitem[Ar]{Ar} S.~Ariki,  \emph{Representation theory of a Hecke algebra of 
$G(r, p, n)$}, J.~Algebra {\bf 177} (1995), 164--185.
%
\bibitem[ArKo]{ArKo} S.~Ariki, K.~Koike, {\em A Hecke algebra of $(\mathbb{Z}/r\mathbb{Z})\wr S_n$ and construction of its irreducible representations}, Adv.~Math.~\textbf{106} (1994), 216--243.
%
\bibitem[Ban]{Ban} E.~Bannai, \emph{Fundamental groups of the spaces of regular orbits of the finite unitary reflection groups of dimension 2},
J.~Math.~Soc.~Japan {\bf 28} (1976), 447--454.
%
\bibitem[Ben]{Ben} M.~Benard, \emph{Schur indices and splitting fields of the unitary reflection
groups}, J.~Algebra {\bf 38} (1976), 318--342.
%
%
\bibitem[Bes1]{Bes1} D.~Bessis, \emph{Sur le corps de d{\'e}finition d'un
groupe de r{\'e}flexions complexe}, Comm.~ Algebra {\bf 25}(8) (1997),
2703--2716.
%
\bibitem[Bes2]{Bes2} D.~Bessis, \emph{Zariski theorems and diagrams for braid groups}, Invent.~Math.~{\bf 145}(3) (2001), 487--507.
%
\bibitem[Bes3]{Bes3} D.~Bessis, \emph{Finite complex reflection arrangements are $K(\pi, 1)$}, Ann.~Math.~{\bf 181}(3) (2015), 809--904.
%
\bibitem[BCCK]{BCCK} C.~Boura, E.~Chavli, M.~Chlouveraki, K.~Karvounis,  \emph{The BMM symmetrising trace conjecture for groups $G_4$, $G_5$, $G_6$, $G_7$, $G_8$}, J.~Symbolic Computation (2019), doi: 10.1016/j.jsc.2019.02.012.
%
\bibitem[Bou]{Bou05} N.~Bourbaki, Lie groups and Lie algebras. Chapters 4--6, Elements of mathematics, English translation of ``Groupes et alg\`ebres de Lie'', Springer, 2005.
%
%
%
\bibitem[BreMa]{BreMa} K.~Bremke, G.~Malle, \emph{Reduced words and a length function for $G(e, 1, n)$}, Indag. Math. {\bf 8} (1997), 453--469.
%
%
\bibitem[BroMa]{BM} M.~Brou\'e, G.~Malle, \emph{Zyklotomische Heckealgebren},  Ast\'erisque  {\bf 212}  (1993), 119--189.
%
\bibitem[BMM]{BMM}  M.~Brou{\'e}, G.~Malle, J.~Michel, \emph{Towards Spetses I}, Trans.~Groups {\bf 4} (1999), 157--218.
%
%
%
\bibitem[BMR]{BMR} M.~Brou{\'e}, G.~Malle, R.~Rouquier, \emph{Complex 
	reflection groups, braid groups, Hecke algebras}, J.~reine angew.
Math. {\bf 500} (1998), 127--190.
%
%
%
\bibitem[Cha1]{Ch17} E.~Chavli, \emph{The BMR freeness conjecture for the tetrahedral and octahedral family}, Comm.~Algebra {\bf 46}(1) (2018), 386--464.
%
\bibitem[Cha2]{Ch18} E.~Chavli, \emph{Universal  deformations  of  the  finite  quotients  of  the  braid  group  on  $3$  strands}, J.~Algebra {\bf 459} (2016), 238--271.
%
%
\bibitem[DMM]{DMM}F.~Digne, I.~Marin, J.~Michel, \emph{The  center  of  pure  complex  braid  groups}, J. Algebra {\bf 347} (2011), 206--213.
%
\bibitem[Fun]{Fun} L.~Funar, \emph{On the quotients of cubic Hecke algebras}, Comm. Math.  Phys. {\bf 173}(3) (1995), 513--558.
%
%
%
\bibitem[LeTa]{lehrer}G.~I.~Lehrer, D.~E.~Taylor,
Unitary reflection groups, volume 20 of
Australian Mathematical Society Lecture Series, Cambridge University Press,
Cambridge, 2009.
%
\bibitem[Mal]{MalleS} G.~Malle, \emph{Spetses}, ICM II, Doc. Math. J. DMV Extra Volume ICM II (1998), 87--96.
%
%
\bibitem[MalMat]{MM98} G.~Malle, A.~Mathas, \emph{Symmetric cyclotomic Hecke algebras}, J. Algebra {\bf 205}(1) (1998), 275--293.
%
%
\bibitem[MalMi]{MM10} G.~Malle, J.~Michel, \emph{Constructing representations of Hecke algebras for complex
reflection groups}, LMS J. Comput. Math. {\bf 13} (2010), 426--450.
%
\bibitem[Mar1]{MarKramm}  I.~Marin, \emph{Krammer representations for complex reflection groups}, J. Algebra {\bf 371} (2012), 175--206.
%
\bibitem[Mar2]{Mar41}  I.~Marin, \emph{The cubic Hecke algebra on at most $5$ strands}, J. Pure Appl. Algebra {\bf 216} (2012), 2754--2782.
%
 \bibitem[Mar3]{Mar43}  I.~Marin, \emph{The freeness conjecture for Hecke algebras of complex reflection groups, and the case of the Hessian group $G_{26}$}, J. Pure Appl. Algebra {\bf 218} (2014), 704--720.
%
\bibitem[Mar4]{MarNew}  I.~Marin, \emph{Proof of the BMR conjecture for $G_{20}$ and $G_{21}$}, J. Symbolic Computation {\bf 92} (2019), 1--14.
%
\bibitem[MarPf]{MaPf}  I.~Marin, G.~Pfeiffer, \emph{The  BMR  freeness  conjecture  for  the  $2$-reflection  groups}, Math.~Comput. {\bf 86} (2017), 2005--2023.
%
\bibitem[MarWa]{Mar46}  I.~Marin, E.~Wagner, \emph{Markov traces on the BMW algebras}, preprint, arXiv:1403.4021.
%
\bibitem[Sage]{sagemath}
    SageMath, The Sage Mathematics Software System (Version 8.1),
   The Sage Developers, 2017, \url{http://www.sagemath.org}.
%
\bibitem[ShTo]{ShTo} G.~C.~Shephard, J.~A.~Todd, \emph{Finite unitary reflection
	groups}, Canad.~J.~Math. {\bf 6} (1954), 274--304.
%
\bibitem[St]{steinberg}R.~Steinberg, \emph{Differential equations invariant under finite reflection groups}, Trans.~ Amer.~ Math.~ Soc. {\bf 112} (1964), 392--400.
%
\bibitem[Tsu]{Tsu} S.~Tsuchioka, \emph{BMR freeness for icosahedral family}, Experimental Mathematics (2018), doi: 10.1080/10586458.2018.1455072.
%
\bibitem[Web]{Web} The project's webpage: \url{https://www.eirinichavli.com/BMMsym13.html}.




\end{thebibliography}
\end{document}